\documentclass[article,ij4uq]{ij4uq}      

\usepackage[hang]{footmisc}
\usepackage[ruled]{algorithm2e}
\usepackage{algorithm2e}
\usepackage{algorithmic}
\fancypagestyle{plain}{%
  \fancyhf{}
  \fancyhead[R]{\small {\it International Journal for Uncertainty Quantification}, x(x): \thepage--\pageref{LastPage} (\myyear\today)}
  \fancyfoot[R]{\small\bf\thepage }
  \fancyfoot[L]{\fottitle}
  }

\renewcommand{\myyear}{2018}
\renewcommand{\today}{}

\usepackage{epstopdf,listings}

\usepackage{tikz}
\usetikzlibrary{positioning}
\usepackage{pgfplots}
\usepackage{pgfplotstable,color}
\usetikzlibrary{patterns}
\usepgfplotslibrary{patchplots}
\definecolor{darkgreen}{rgb}{0,0.7,0}
\definecolor{lightblue}{rgb}{0.8,0.8,1.0}
\usetikzlibrary{decorations.text}
\usetikzlibrary{decorations.pathmorphing}

\pgfplotscreateplotcyclelist{mycolor}{%
blue,every mark/.append style={fill=blue!80!black},mark=*\\%
red,every mark/.append style={fill=red!80!black},mark=square*\\%
darkgreen,every mark/.append style={fill=darkgreen!80!black},mark=otimes*\\%
black,mark=star\\%
blue,every mark/.append style={fill=blue!80!black},mark=diamond*\\%
red,densely dashed,every mark/.append style={solid,fill=red!80!black},mark=*\\%
darkgreen,densely dashed,every mark/.append style={
solid,fill=darkgreen!80!black},mark=square*\\%
black,densely dashed,every mark/.append style={solid,fill=gray},mark=otimes*\\%
blue,densely dashed,mark=star,every mark/.append style=solid\\%
red,densely dashed,every mark/.append style={solid,fill=red!80!black},mark=diamond*\\%
}

\pgfplotsset{cycle list name={mycolor}}
\pgfplotsset{
/pgfplots/bar cycle list/.style={/pgfplots/cycle list={%
{blue,fill=blue,mark=none},%
{red,fill=red,mark=none},%
{darkgreen,fill=darkgreen,mark=none},%
{black,fill=black,mark=none},%
}
},
}

\newcommand{\bnd}{\Gamma}

\newcommand{\cov}[1]{\mbox{Cov} \left[#1 \right]}

\newcommand{\E}[1]{\mathbb{E}\left[#1 \right]}
\newcommand{\R}{\mathbb{R}}

\renewcommand{\vec}[1]{{\bf #1}}
\newcommand{\tens}[1]{{\bf #1}}
\newcommand{\dom}{\mathcal{D}}

\newcommand{\kernel}{k}

\newcommand{\pdim}{D}

\newcommand{\KL}{Karhunen-Lo\`eve\ }

\newcommand{\paradom}{\mathcal{Y}}
\newcommand{\Nsamp}{N}
\newcommand{\Ngrid}{M}
\newcommand{\kernelscale}{\zeta}

\makeatletter
\newcommand{\vast}{\bBigg@{3}}
\newcommand{\Vast}{\bBigg@{4}}
\makeatother

\DeclareMathOperator{\spn}{span}

\DeclareMathOperator\erf{erf}

\newcommand{\reviewerone}[1]{{#1}}
\newcommand{\reviewertwo}[1]{{#1}}
\newcommand{\crcom}[1]{{#1}}
\begin{document}

\volume{Volume x, Issue x, \myyear\today}
\title{Kernel--based stochastic collocation for the random two--phase Navier-Stokes equations}
\titlehead{Kernel--based UQ for two--phase Navier Stokes}
\authorhead{M.~Griebel, C.~Rieger, \& P.~Zaspel}
\author[1]{M.~Griebel}
\author[2]{C.~Rieger}
\corrauthor[3]{P.~Zaspel}
\corremail{peter.zaspel@unibas.ch}
\corraddress{Department of Mathematics und Computer Science, University of Basel, Spiegelgasse 1, 4051 Basel, Switzerland}
\address[1]{Institute for Numerical Simulation, Bonn University, \crcom{Endenicher Allee 19b}, D-53115 Bonn, Germany \& 
Fraunhofer Institute for Algorithms and Scientific Computing SCAI, Schloss Birlinghoven, D-53754 Sankt Augustin, Germany
}
\address[2]{Institute for Numerical Simulation, Bonn University, \crcom{Endenicher Allee 19b}, D-53115 Bonn, Germany \&  
Department of Mathematics, RWTH Aachen University, Schinkelstr. 2, D-52062 Aachen, Germany
}
\address[3]{Department of Mathematics and Computer Science, University of Basel, Spiegelgasse 1, 4051 Basel, Switzerland}

\dataO{mm/dd/yyyy}
\dataF{mm/dd/yyyy}

\abstract{In this work, we apply stochastic collocation methods  with radial kernel basis functions for an uncertainty quantification of the random incompressible two-phase Navier--Stokes equations. Our approach is non-intrusive and we use the existing fluid dynamics solver NaSt3DGPF to solve the incompressible two-phase Navier--Stokes equation for each given realization. We are able to empirically show that the resulting kernel-based stochastic collocation is highly competitive in this setting and even outperforms some other standard methods.  
}

\keywords{ stochastic collocation, incompressible two--phase Navier--Stokes, uncertainty quantification}
\maketitle


\section{Introduction}
In this paper, we apply uncertainty quantification to the
large-scale complex fluid dynamics problem of incompressible two-phase flows modeled by the three-dimensional random two-phase Navier--Stokes equations.
The two-phase Navier--Stokes equations describe the interaction of two non-mixing fluids like water
and oil or water and air (at low Mach numbers). It has important applications ranging from fluvial construction analysis to flows in chemical bubble reactors. In fluvial construction analysis, a quantification of uncertainties is crucial for public safety. In chemical bubble column reactors, we can use uncertainty quantification to obtain a stochastic homogenization of the perturbation of a liquid in presence of many rising bubbles. This is important for large-scale chemical process optimization.

Depending on the respective applications, we treat densities, viscosities or volume forces as stochastic values or fields. This renders also the systems' solution, i.e.~the velocity field, the pressure field and the liquid-liquid interface as stochastic quantities. Moreover, quantities of interest computed from these solution fields thereby also become stochastic. After the numerical solution of the overall problem, its uncertainty quantification can be achieved by a stochastic moment analysis, which will include the evaluation of the e.g., first stochastic moment. 
Methods for stochastic moment analysis are intrusive (e.g.~stochastic Galerkin \cite{Ghanem1991,
Schwab2011}) or 
non-intrusive (e.g.~Monte Carlo, quasi-Monte Carlo \cite{Graham2011}, multi-level Monte Carlo \cite{Barth2011} and (generalized) polynomial chaos \cite{Jakeman2008} and stochastic collocation
\cite{Babuska2010}). 
Moreover, there is some related work for stochastic moment analysis in computational fluid dynamics applications, with examples in
groundwater flow \cite{Ganis2008}, incompressible flows
\cite{
LeMaitre2010,
Schicktoappear}. 
However, to our knowledge, the application of the two-phase incompressible Navier--Stokes equations has never been considered
before.

Our approach to solve the random two-phase Navier--Stokes equations will be based on a \textit{non-intrusive stochastic collocation} approach. This enables us to re-use our existing two-phase Navier--Stokes solver NaSt3DGPF \cite{Croce2009,Dornseifer1998}. 
This code covers applications such as river simulations in presence of hydraulic constructions and, more recently, sediment transport
 \cite{Burkow2013}
and non-Newtonian flows
\cite{Griebel2013}. 
It has also been parallelized on CPU and GPU clusters \cite{Zaspel2013}.

In stochastic collocation, spectral (sparse) tensor-product approximation \cite{Babuska2010} is widely used. However, since regularity results are in general not available for the two--phase Navier--Stokes equations, we can in general not expect to have the high parametric regularity required by some of the recent spectral approximation approaches (even in the regime of moderately high Reynolds numbers). Moreover, the number of samples in grid-based approximation approaches, even for sparse grids, can be prohibitively high due to its dyadic construction. This is problematic for our application, since already one single deterministic high-resolution two--phase Navier--Stokes simulation often requires computational resources in the range of hours to days on a parallel computer. 
Therefore, we introduce a meshfree approach for the treatment of the stochastic variables, namely the radial basis function (RBF)
kernel-based stochastic collocation, to achieve both high asymptotic convergence rates and good pre-asymptotic behavior for
stochastic moment analysis.
We note the close relationship of kernel-based
approximation to kriging \reviewertwo{\cite{Krige1951,Beers2004}} with its
low error for few collocation points and to Gaussian process regression \reviewertwo{\cite{Rasmussen2005}} with its profound stochastic framework. 
Related recent work on kernel-based collocation covers the
approximate solution of stochastic partial differential equations
\cite{Cialenco2012}, the special case of
an elliptic random PDE \cite{Fasshauer2013} by an intrusive method, the case of parametric partial differential equations \cite{Griebel2017} and the parallel treatment of large random partial differential equations, see \cite{Zaspel2015}.

In our numerical results with several rising bubble test cases, we focus on a small number of random parameters, which are of nearly equal importance to the simulation. 
This is justified by an engineering perspective where a few physical parameters are considered with small fluctuations.
A small number of parameters is also justified by practical limitations since we would never be able to sample enough solutions (or store them) to sufficiently resolve a larger and higher dimensional parameter domain. 
As discussed before, we do not assume or expect to have high spatial/temporal regularity for the quantities of interest. 
Therefore, we employ algebraically smooth kernels to this situation. Algebraically smooth kernels are known to achieve a convergence rate which corresponds to the minimum of the smoothness of the kernel and the smoothness of the function to be approximated, see \cite{NarcowichWard04}. Hence, by using a kernel approach with higher algebraic regularity, we employ a method which is capable to exploit the unknown smoothness of the parametric function which needs to be reconstructed. 
Moreover, in the isotropic case where all dimensions in the parameter space are nearly equally important, we know that kernel methods are quasi-optimal, see \cite{Wendland2004}. From a practical point of view, kernel methods are easy to implement compared to more sophisticated approximation methods such as multilevel or sparse grid constructions. Finally, kernel methods allow to increase the number of sampling points by just any arbitrary number and we do not need to stick to prime numbers (as in some QMC methods), grid sizes (as in some sparse grids) or dimensions of polynomial spaces.  
In addition, we use the Gaussian kernel if the quantity of interest might depend smoothly on the parameters as we expect it for integrated quantities of interest. In such a situation, kernel methods are able to provide exponential convergence rates from the very beginning. We confirm these statements by numerical results.
In particular, we show empirically that kernel-based stochastic collocation methods
allow to outperform the algebraic convergence rates of isotropic (quasi-)Monte Carlo techniques.
Moreover, in a direct comparison to a (sparse) spectral tensor-product method, in which we consider an integrated quantity of interest, the kernel-based method shows exponential convergence rates. Hence, our approach is able to deliver a decent approximation with very few deterministic solutions of the two-phase Navier--Stokes equations.

\newpage
The remainder of this work is organized as follows. 
Section~\ref{sec:modelNavierStokes} introduces the random two-phase Navier--Stokes equations including details on their numerical treatment.
Section~\ref{chap:RBFKernelMethods} discusses the RBF kernel-based stochastic collocation. 
The main results, i.e., empirical convergence studies, are presented in Section~\ref{sec:numericalResults}. Finally, Section~\ref{sec:conclusions} gives some conclusions and a short outlook.

\section{Random two-phase Navier--Stokes equations}
\label{sec:modelNavierStokes}
The application problem motivating this work is a random version of the two-phase incompressible Navier--Stokes equations.
They model the interaction of two incompressible fluids which do not mix but
remain disjoint with a common interface. Classical examples for such fluid-fluid
systems in real world are oil and water or water and air, noting that it is
usual to assume air to be incompressible for the discussed test
cases. 
The random model discussed here is derived from the deterministic two-phase Navier--Stokes equations as discussed e.g.~in~\cite{Croce2009}. Further references are
\cite{Sussman1994,Gross2011}. 

Formally, the introduction of the random two-phase Navier--Stokes equations would start from a probability space (possibly infinite-dimensional) with stochastic
parameters in which all stochastic variations are expressed. However, 
the finite-noise approximation leads us to consider only a $\pdim$-dimensional parameter space $\paradom \subset \R^{\pdim}$, cf.~e.g.~\cite{Babuska2010}. The parameter space $\paradom$ is considered to be a measure space $(\paradom, \mathcal{B},\rho\,d\vec{y})$ with density $\rho$ and the usual Borel $\sigma$--algebra $\mathcal{B}$. Here we directly introduce the
\textit{finite-dimensional} (parametric) model of the random two-phase Navier--Stokes equations, noting
that the parameter space $\paradom\subset \R^{\pdim}$ might be still high-dimensional, i.e., $\pdim$ might be large.

\subsection{Finite-dimensional random model}
\begin{figure}
\begin{center}
\begin{tikzpicture}
 \filldraw[lightblue] (0,1.5) sin (0.5,1.75) cos (1,1.5) sin (1.5,1.25) cos
  (2,1.5) sin (2.5,1.75) cos (3,1.5) sin (3.5,1.25) cos (4,1.5)
              sin (4.5,1.75) cos (5,1.5) sin (5.5,1.25) cos (6,1.5)
              sin (6.5,1.75) cos (7,1.5) -- (7,0) -- (0,0) -- (0,1.5);
  \draw[lightblue] (3.5,1) node[anchor=north]
  {\textcolor{black}{$\dom_1$}}; \draw (3.5,3) node[anchor=north]
  {$\dom_2$};
  \draw (0,0) rectangle (7,4);
  \draw[gray] (7,2.6) -- (6.5,3.2) node[left] {\textcolor{black}{$\bnd$}};
  \draw[gray] (2.5,1.75) -- (1.75,2) node[above] {\textcolor{black}{$\bnd_f$}}; 
\end{tikzpicture}
\end{center}\vspace*{-0.75em}
\caption{\label{sec:ModelNavierStokesDomains}The domain $\dom$ with boundary
$\bnd$ is subdivided into two distinct fluid phase domains $\dom_1$ and
$\dom_2$ and the fluid-fluid interface $\bnd_f$ in the two-phase Navier--Stokes
equations.}
\end{figure}
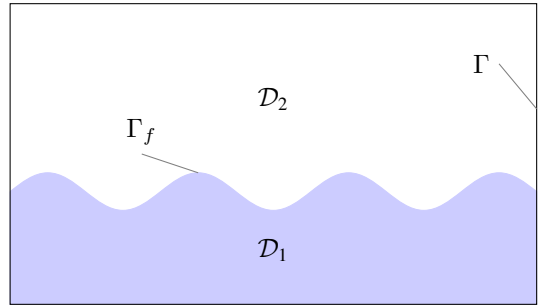

Now, we present the finite-dimensional random two-phase Navier--Stokes
equations with respect to the parameter space $(\paradom, \mathcal{B},\rho\,d\vec{y})$. The physical space
$\dom\subset\R^3$ is a connected domain with boundary $\bnd=\partial\dom$. Two sub-domains
$\dom_1,\dom_2$ identify the two fluid phases. Technically, they
are time- and parameter-dependent, i.e, $\dom_{i}=\dom_{i}(\vec{y},t)$ for $i=1,2$
with $\dom_1(\vec{y},t)\cap\dom_2(\vec{y},t)=\emptyset$ for all $t\in [0,T]$ and $T\in\R_+$ the final
simulation time. The full domain $\dom$ is covered by $\dom_1$,
$\dom_2$ and the time- and parameter-dependent fluid-fluid separation interface
$\bnd_f(\vec{y},t)$. Thus there holds $\dom = \dom_1(\vec{y},t) \cup \dom_2(\vec{y},t) \cup \bnd_f(\vec{y},t)$,
cf.~Figure~\ref{sec:ModelNavierStokesDomains}. Moreover, the outer domain boundary is decomposed into $\bnd_{1}:=\bnd \cap \partial \dom_1$ and $\bnd_{2}:=\bnd \cap \partial \dom_2$. We note that by this construction the whole domain $\dom$ is fixed but the two sub-domains and the free surface separating them are parameter-dependent.
In each of the two sub-domains $\dom_{i}$, $i=1,2$, the finite-dimensional random system of the
two-phase Navier--Stokes equations reads as
\begin{linenomath*}
\begin{eqnarray}	
\label{eq:ModelDeterministicRandomNSmomentum}\rho_i(\vec{y})\frac{D\vec{u}_i}{D
t} =&\hspace*{-0.75em} \nabla\cdot \mu_i(\vec{y})\tens{S}_i - \nabla p_i +
\rho_i(\vec{y})\vec{g}(\vec{y}) & \mbox{in}\
\paradom\times\dom_i(\vec{y})\times [0,T],\\
\label{eq:ModelDeterministicRandomNScontinuity}\nabla\cdot\vec{u}_i =& 0 & \mbox{in}\
\paradom\times\dom_i(\vec{y})\times [0,T],\\
\vec{u}_i = & \vec{u}_{0_i}(\vec{y},\vec{x}) & \mbox{in}\
\paradom\times\dom_i(\vec{y})\times \{0\},\\
\label{eq:ModelDeterministicRandomNSbc1}\mathcal{B}\vec{u}_i =&
\vec{b}_\bnd & \mbox{on}\ \paradom\times\bnd\times [0,T],\\
\label{eq:ModelDeterministicRandomNSbc2} \frac{\partial p_i}{\partial\vec{n}_\bnd} = & 0 &
\mbox{on}\ \paradom\times\bnd\times [0,T]\,,\\
\label{eq:ModelDeterministicRandomNSinterfaceVelocity}\vec{u}_1 =& \vec{u}_2 & \mbox{on}\
\paradom\times\bnd_f(\vec{y})\times [0,T],\\
\label{eq:ModelDeterministicRandomNCinterfaceJump}\reviewertwo{\left[ {\bf
T \vec{n}}\right]} =& \sigma\kappa\vec{n} &
\mbox{on}\ \paradom\times\bnd_f(\vec{y})\times [0,T].
\end{eqnarray}
\end{linenomath*}

The main part of the random two-phase Navier--Stokes equations are the
momentum equation \eqref{eq:ModelDeterministicRandomNSmomentum} and the continuity
equation \eqref{eq:ModelDeterministicRandomNScontinuity}. \reviewertwo{The first equation is decomposed into the material derivative 
$\frac{D\vec{u}_i}{Dt} := \partial_t \vec{u}_i + (\vec{u}_i\cdot \nabla)\vec{u}_i$ and terms involving the viscosity (with $\tens{S}_i := \nabla\vec{u}_i + \{\nabla\vec{u}_i\}^T$), pressure and volume forces. The second equation represents the incompressibility
constraint for both fluids.}

Both fluids interact with respect to a random volume force $\vec{g}:\paradom\times\dom\rightarrow\R^3$, e.g.~gravity with some perturbation.\reviewertwo{\footnote{\reviewertwo{While a random volume force might not have an immediate physical correspondence, the choice of a random right-hand side in literature on uncertainty quantification for model partial differential equations is very common. Therefore, this case is addressed here, too.}}}
At the fluid-fluid
interface $\bnd_f$, there holds the jump condition
\eqref{eq:ModelDeterministicRandomNCinterfaceJump} for the \reviewertwo{stress tensor $\tens{T}_i :=
-p_i\tens{I} + \mu_i\tens{S}_i$,\linebreak $\tens{T}_i:\paradom\times\dom_i\times [0,T]\rightarrow\R^{3\times 3}$}, with
\reviewertwo{$\left[ \tens{T}\vec{n}\right]$ being the jump $(\tens{T}_1\vec{n}-\tens{T}_2\vec{n})$} across the
interface \reviewertwo{and $\vec{I}\in\R^{3\times 3}$ being the identity matrix}. The continuity condition for the velocitities
$\vec{u}_1$ and $\vec{u}_2$ at the interface is given in
\eqref{eq:ModelDeterministicRandomNSinterfaceVelocity}. Finally, $\sigma\in\R$ is the surface
tension coefficient, $\kappa:\paradom\times\bnd_f(\vec{y})\times [0,T]\rightarrow\R$ is the curvature of $\bnd_f$ and
$\vec{n}:\paradom\times\bnd_f(\vec{y})\times [0,T]\rightarrow\R^3$ is the surface normal of the interface.
\reviewertwo{For simplicity, we keep $\sigma$ to be deterministic.}

This system is augmented by proper initial conditions for the
velocity field by $\vec{u}_{0_i}:\dom_i\rightarrow\R^3$\reviewertwo{, by boundary
conditions \eqref{eq:ModelDeterministicRandomNSbc1} for the velocities and by approximate boundary conditions \eqref{eq:ModelDeterministicRandomNSbc2} for the pressures}. 
Velocity boundary conditions are here denoted for the sake of simplicity by some general
boundary operator $\mathcal{B}$ and the space-time-dependent right-hand side function
$\vec{b}_\bnd$.  
For all parameters $\vec{y}\in\paradom$, equations \eqref{eq:ModelDeterministicRandomNSmomentum}--\eqref{eq:ModelDeterministicRandomNCinterfaceJump} are solved for the velocity fields $\vec{u}_i:\paradom\times\dom_i\times
[0,T]\rightarrow\R^3\,\,[m/s]$ and pressures $p_i:\paradom\times\dom_i\times
[0,T]\rightarrow\R\,\,[kg/(m\cdot s^2)]$. 
The two important material properties for incompressible fluids are the
(random) subdomain-wise constant densities $\rho_i:\paradom\rightarrow\R\,\,[kg/m^3]$ and viscosities
$\mu_i:\paradom\rightarrow\R\,\,[kg/(m\cdot s)]$.
Note that it is common to have no initial
conditions for the pressure, since the pressure is usually understood as a Lagrange multiplier. Moreover, the solution method
applied in this work does not require initial conditions for pressure.

Altogether, this formulation of the random model explicitly introduces a parametric/random dependence of the two-phase Navier--Stokes equations in the densities, viscosities, the volume force and the initial condition. This way, all other quantities become dependent on the random/parametric input.

\subsection{Continuous formulation using random level-sets}
We now follow the common approach to
introduce a level-set function $\phi$ to distinguish the two fluid phases. For further details, see \cite{Croce2009}. This way, the equations can be formulated with respect to the common domain $\dom$. 
Deterministic level-set techniques have been proven useful in the context of two-phase flows is \cite{Sussman1994}. In our situation the level set function $\phi$ is now random/parametric. To this end, $\phi:\paradom\times\dom\times [0,T]\rightarrow\R$ is a signed
distance function implicitly defining the two random domains $\dom_1$, $\dom_2$ as
\begin{linenomath*}
\begin{equation}\label{levelset}
\phi(\vec{y},\vec{x},t) \left\{
\begin{array}{ll}
        < 0 & \:\:\mbox{if}\,\,\, \vec{x}\in\dom_1(\vec{y},t),\\
    = 0 & \:\:\mbox{if}\,\,\, \vec{x}\in\bnd_f(\vec{y},t),\\
    > 0 & \:\:\mbox{if}\,\,\, \vec{x}\in\dom_2(\vec{y},t).
\end{array}\right.
\end{equation}
\end{linenomath*}
Also, $\phi$ obeys the Eikonal equation 
\begin{linenomath*}
\reviewerone{\begin{equation}
|\nabla_{\vec{x}}\phi(\vec{y},\vec{x},t)| = 1 \quad \text{for almost all }(\vec{y},\vec{x},t) \in \paradom\times\dom\times [0,T],
\end{equation}}
\end{linenomath*}
which makes it for fixed parameter $\vec{y} \in \paradom$ a distance function. The free surface $\bnd_f$, i.e., the (random)
interface between both fluids, is given by
\begin{linenomath*}
\reviewerone{\begin{equation}\bnd_f (\vec{y},t) = \{\vec{x} :
\phi(\vec{y},\vec{x},t) = 0\}\,.\end{equation}}
\end{linenomath*}
Based on the Continuum Surface Force scheme \cite{Brackbill1992} it is possible
to reformulate the discontinuous equation system
\eqref{eq:ModelDeterministicRandomNSmomentum}--\eqref{eq:ModelDeterministicRandomNCinterfaceJump} into a continuous
representation \cite{Sussman1994}. This works analogously for the random/parametric case. 
In the remainder of this section, we understand all \reviewertwo{equations} in the sense that they hold for (almost) every parameter value $\vec{y} \in \paradom$.
We here follow \cite{Croce2009} and introduce, by slightly abusing notation, domain-dependent densities and
viscosities
\begin{linenomath*}
\reviewerone{\begin{equation}\rho_\phi(\vec{y}) := \rho_1(\vec{y}) + (\rho_2(\vec{y}) - \rho_1(\vec{y}) )\, H (\phi(\vec{y})),\end{equation}}
\end{linenomath*}
\begin{linenomath*}
\reviewerone{\begin{equation}\mu_\phi(\vec{y}) := \mu_1(\vec{y}) + (\mu_2(\vec{y}) - \mu_1(\vec{y}) )\, H (\phi(\vec{y}))\end{equation}}
\end{linenomath*}
with the Heaviside step function
\begin{linenomath*}
\reviewerone{\begin{equation}
H(\phi ) := \left\{
\begin{array}{ll}
        0 & \:\:\mbox{if}\,\,\, \phi < 0,\\
    \frac{1}{2} & \:\:\mbox{if}\,\,\, \phi = 0,\\
    1 & \:\:\mbox{if}\,\,\, \phi > 0. 
\end{array}
\right. 
\end{equation}}
\end{linenomath*}
This function is smoothed out in an $\epsilon$-environment of the free
surface leading to jump-free functions $H^\epsilon$, $\rho_{\phi}^\epsilon$ and
$\mu^\epsilon(\phi)$, cf.~\cite{Croce2009} for more details.
It is then possible to derive the initial-boundary value problem
\begin{linenomath*}
\begin{eqnarray}
\label{eq:NSDeterministicRandomContinuousMomentum}\rho_\phi^\epsilon \frac{D \vec{u}}{D t}  =&
\hspace*{-0.75em}\nabla\cdot (\mu_\phi^\epsilon \tens{S}) - \nabla p -
\sigma\kappa_\phi\delta^\epsilon \nabla\phi + \rho_\phi^\epsilon
\vec{g}\hspace*{-0.45em}  & \mbox{in}\ \paradom\times\dom\times [0,T],\\
\label{eq:NSDeterministicRandomContinuousContinuity}\nabla\cdot\vec{u} =&  0 & \mbox{in}\
\paradom\times\dom\times [0,T],\\
\label{eq:NSDeterministicRandomLevelSetTransport}\hspace*{-0.75em}\partial_t\phi + \vec{u}\cdot\nabla\phi =& 0  &
\mbox{in}\ \paradom\times\dom\times [0,T],\\
\label{eq:NSDeterministicRandomContinuousEikonal}\left|\nabla\phi\right| =& 1& \mbox{on}\
\paradom\times\bnd\times [0,T],\\
\label{eq:NSDeterministicRandomContinuousinitial}\vec{u} =& \vec{u}_0 & \mbox{in}\ \paradom\times\dom\times\{0\},\\
\mathcal{B}\vec{u} =& \vec{b}_\bnd & \mbox{on}\
\paradom\times\bnd\times [0,T],\\
\frac{\partial p}{\partial\vec{n}_\bnd} =&
0 & \mbox{on}\ \paradom\times\bnd\times [0,T],\\
\phi =& \phi_0 & \mbox{in}\ \paradom\times\dom\times
\{0\},\\
\label{eq:NSDeterministicRandomContinuousLast}\frac{\partial \phi}{\partial
\vec{n}_\bnd} = & 0 & \mbox{on}\ \paradom\times\bnd\times[0,T],
\end{eqnarray}
\end{linenomath*}
with $\vec{u}:\paradom\times\dom\times [0,T]\rightarrow \R^3$ and $p:\paradom\times\dom\times
[0,T]\rightarrow \R$ the velocity and pressure fields, respectively, defined on the full domain
$\dom$. It turns out that the jump condition for the stress tensor translates to the volume force
$-\sigma\kappa_\phi^\epsilon \delta_\phi^\epsilon \nabla\phi$, where
$\delta_\phi^\epsilon$ denotes a smoothed out Dirac functional and $\kappa_\phi$
is the curvature of the free surface. It is given in the level-set case (cf.
\cite{Osher1988}) as 
\begin{linenomath*}
\reviewerone{\begin{align}\kappa_\phi(\vec{y},\vec{x},t) =
\nabla_{\vec{x}}\cdot\frac{\nabla_{\vec{x}}\phi(\vec{y},\vec{x},t)}{\|\nabla_{\vec{x}}\phi(\vec{y},\vec{x},t)\|}  \quad \text{for } \vec{x} \in
\bnd_f(\vec{y},t) \,.\end{align}}
\end{linenomath*}
\reviewertwo{The transport equation \eqref{eq:NSDeterministicRandomLevelSetTransport} governs the evolution of the interface.}
Equation \eqref{eq:NSDeterministicRandomContinuousEikonal} describes the re-initialization of the level-set function and introduces \reviewertwo{a nonlinearity}. It will be treated by an iterative method, see \cite{Croce2009} for further details.
Again, the system is augmented by proper initial and boundary conditions, see \eqref{eq:NSDeterministicRandomContinuousinitial}--\eqref{eq:NSDeterministicRandomContinuousLast}.
Altogether, the equations \eqref{eq:NSDeterministicRandomContinuousMomentum}--\eqref{eq:NSDeterministicRandomContinuousLast} describe a Continuum Surface Force formulation of the parametric Navier--Stokes equation. As already mentioned, this system of equations holds (in the sense of almost everywhere) for each parameter. For a fixed parameter, the system of equations reduces to a common deterministic two--phase Navier--Stokes equation.

\subsection{Numerical treatment of the two-phase Navier--Stokes equations}
The stochastic collocation approach only requires us to consider the solution of the deterministic two-phase Navier--Stokes equations for a set of fixed parameter values $\vec{y}\in\paradom$. To this end, classical solution techniques can be used. 
In practice, we employ the flow solver NaSt3DGPF \cite{Dornseifer1998,Croce2009}, which implements the pressure correction approach \cite{Chorin1968} to solve the deterministic two-phase Navier--Stokes equations over time. Details can be found in \cite{Croce2009}.
The main features of the code are as follows:

Finite differences/volumes are used as discretization method in space with a
classical marker-and-cell (MAC) staggered uniform grid.
The velocity components are therefore discretized on the centers of the cell faces
whereas pressure and level-set function are discretized on the centers of the
cells.
Wherever it is necessary to evaluate quantities e.g.~from cell centers on
the cell faces, higher-order interpolation is used. The convective terms
of the momentum equations
are discretized by the fifth-order weighted essentially non-oscillatory (WENO) scheme. 
The diffusion term is computed by second-order central differences. The WENO scheme
is also applied to the gradient evaluation in the reinitialization equation and to the
transport term in the level-set advection.
The pressure Poisson equation is discretized with a standard seven-point second
order stencil and solved with a Jacobi-preconditioned conjugate gradient (CG)
method. 
As time integrator, a second order  
Adams--Bashforth scheme with an adaptive time
step selection mechanism is employed, cf.~\cite{Dornseifer1998}, which obeys the Courant--Friedrichs--Lewy
(CFL) condition.

As a consequence of this discretization and the numerical solution, note that we have at no point of our stochastic collocation procedure access to the true solutions $\vec{u}$, $p$ and $\phi$. Instead, we will restrict ourselves to the numerical approximations $\vec{u}^{h}$, $p^{h}$ and $\phi^{h}$ where $h$ indicates spatial and temporal discretization.

\section{Kernel-based stochastic collocation for stochastic moment analysis
}
\label{chap:RBFKernelMethods} 

\subsection{Objective}
The numerical solution to the random two--phase Navier--Stokes system with random level-sets consists of three components, namely a time--dependent velocity random vector-field $\vec{u}^{h}:\paradom\times \bar{\dom}\times [0,T]\rightarrow\R^3$ 
with $\vec{u}^{h}=(u^{h}_1, u^{h}_2, u^{h}_3)^T$, a time--dependent pressure random field
$p^{h}:\paradom\times \bar{\dom}\times [0,T]\rightarrow \R$
and time--dependent random domains $\dom^{h}_{i}(\vec{y}, t)$ for $i=1,2$, which we represent by the time--dependent level-set random field 
$
\phi^{h}:\paradom\times \bar{\dom}\times [0,T]\rightarrow \R
$.

It is common not to consider the full solution but some \emph{quantities of interest}, which are to be derived from these numerical solutions. To this end, we consider in the following two quantities of interest. The first one is just the first component of the solution velocity field
\begin{linenomath*}\vspace*{-1em}
\reviewerone{\begin{equation}u^{h}_1:\paradom\times \bar{\dom}\times [0,T]\rightarrow\R\,.\vspace*{-1em}\end{equation}}
\end{linenomath*}
This is a parameter-dependent real-valued field discretized in space and time. The second quantity of interest is the center of mass 
\begin{linenomath*}
\reviewerone{\begin{equation}\vec{c}^{h}: \paradom\times [0,T]\rightarrow \R^3\end{equation}}
\end{linenomath*}
with $\vec{c}^{h}:=(c^{h}_1, c^{h}_2, c^{h}_3)^\top$ of the second fluid phase and
\begin{linenomath*}
\reviewerone{\begin{align}
\vec{c}^{h}(\vec{y},t) := \frac{1}{\mbox{Vol}(\dom^{h}_2(\vec{y},t))} \int_{\dom^{h}_2(\vec{y},t)} \vec{x}\,d\vec{x}=\left(\int_{\dom^{h}_2(\vec{y},t)} 1\,d\vec{x} \right)^{-1}\int_{\dom^{h}_2(\vec{y},t)} \vec{x}\,d\vec{x}\ .
\end{align}}
\end{linenomath*}
Obviously, this is a vector-valued parameter-dependent quantity in time. We will here only focus on the second component, namely $c^{h}_2$.

For both quantities of interest, we are interested to compute the first statistical moment, that is we only compute  
\begin{linenomath*}
\begin{equation}\label{exqoi}
	\E{u^{h}_1}(\vec{x},t) = \int_{\paradom}u^{h}_1(\vec{y},\vec{x},t)\rho(\vec{y})\,d\vec{y}\vspace*{-1em}
\end{equation}
\end{linenomath*}
and
\begin{linenomath*}\vspace*{-1em}
\begin{align}\label{eq:firstMomentCenter}
\E{c^{h}_2}(t) &= \int_{\paradom}c^{h}_2(\vec{y},t)\rho(\vec{y})\,d\vec{y}\,.\end{align}
\end{linenomath*}
While all quantities of interest and their statistical moments are here given as time-dependent values, we will in practice always evaluate them for the final time $t=T$ of a given Navier--Stokes simulation only. This gives us a space--discrete function, cf.~Eqs.~\eqref{exqoi}, or a single number, cf.~Eq.~\eqref{eq:firstMomentCenter}.

\subsection{Stochastic collocation}
\label{sec:stochCollocLagrange}
\reviewerone{To underline the generality of the kernel-based stochastic collocation, we will consider a general quantity of interest $u:\paradom\times\bar{\dom}\times [0,T]\rightarrow\R$ for the rest of this section. We leave it to the reader to associate its approximation with the concrete examples of $u_1^h$ and $c_2^h$, which were discussed before.}

Following \cite{Babuska2010}, stochastic collocation evaluates the quantities of interest with respect to $\vec{y}$ in 
collocation points
\begin{linenomath*}
\begin{equation}\label{collocpoints}
Y_{\Nsamp}:=\left\{\vec{y}_1,\ldots,\vec{y}_{\Nsamp}\right\} \subset \paradom\,.
\end{equation}
\end{linenomath*}
\reviewerone{The point evaluations $u(\vec{y}_i,\vec{x},t)$ are used
to numerically approximate the continuous function}
\begin{linenomath*}
\begin{align}\label{eq:lagrangeInterpolationForStochasticCollocation}
u(\vec{y},\vec{x},t) &\approx \sum_{i=1}^{\Nsamp}
u(\vec{y}_i,\vec{x},t) L_i(\vec{y}) & \forall
\vec{y}\in\paradom\,.\end{align}
\end{linenomath*}
This is done with a Lagrange basis in a certain function (approximation) space $\mathcal{H}$ with respect to $\paradom$, where the space $\mathcal{H}$ still needs to be specified, i.e.,
\begin{linenomath*}
\reviewerone{\begin{equation}
\left\{L_i\right\}_{i=1}^{\Nsamp},\quad L_i\in \mathcal{H},\quad
L_i:\paradom\rightarrow\R,\quad \mbox{with}\ L_i(\vec{y}_j) =
\left\{\begin{array}{cc} 1 & i=j,\\0 & i\neq j.\end{array}\right.
\end{equation}}
\end{linenomath*}

Since each evaluation of the quantities of interest for a fixed parameter value $\vec{y}\in Y_{\Nsamp}$ involves \reviewerone{e.g.}~a whole deterministic fluid-dynamics simulation, we do not want the set $Y_{\Nsamp}$ to be of any regular grid-like structure. Grid like structures suffer in particular in a mesh refinement step. Typically, not just one point can be added at a time but several points have to be added in order to maintain the grid structure. Instead, we aim for a meshfree collocation procedure.

Another issue is that we do not know the smoothness of the dependence of the quantities of interest on the parameters. 
This makes the use of approximation by reproducing kernels and in particular radial basis functions favorable since these methods are known to adapt to the smoothness of the function which is reconstructed once the radial basis function is chosen smooth enough and the point set $Y_{\Nsamp}$ is not too far from quasi-uniformity, see \cite{NarcowichWard04}. Consequently, we will choose the function space $\mathcal{H}$ to be a reproducing kernel Hilbert space (RKHS).

\subsection{Reproducing kernel Hilbert spaces and native spaces}
\label{sec:RKHS} 
Following \cite{Wendland2004}, a reproducing kernel Hilbert space is defined as a Hilbert space of functions
\begin{linenomath*}
\reviewerone{\begin{equation}
\mathcal{H}_{\kernel}(\paradom)\subseteq \left\{ f: \paradom\rightarrow\R\, \middle|\,\emptyset\neq\paradom\subseteq\R^{\pdim}\right\}\,
\end{equation}}
\end{linenomath*}
which posseses a kernel function $\kernel:\paradom\times\paradom\rightarrow\R$ that satisfies\vspace*{-1em} 
\begin{linenomath*}
\reviewerone{\begin{align}
\kernel (\cdot,\vec{y})&\in\mathcal{H}_{\kernel}(\paradom) \quad \text{for all } \vec{y}\in\paradom \quad \text{and}\\
f(\vec{y}) &= \left(f, \kernel (\cdot,\vec{y})\right)_{\mathcal{H}_{\kernel}(\paradom)} \quad \text{for all } f\in\mathcal{H}_{\kernel}(\paradom) \  \text{and all }  \ \vec{y}\in\paradom.\vspace*{-1em}
\end{align}}
\end{linenomath*}
This kernel $\kernel$ is called the reproducing kernel.
\reviewertwo{For two functions $f,g\in \mathcal{H}_{\kernel}$ with $f = \sum_{j=1}^N \alpha_{j}^{(f)} \kernel(\cdot,\vec{y}_j)$ the inner product is defined as $(f, g)_{\mathcal{H}_{\kernel}} := \sum_{j=1}^{N} \sum_{j^\prime=1}^N \alpha_{j}^{(f)} \alpha_{j^\prime}^{(g)} \kernel(\vec{y}_j,\vec{y}_{j^\prime})$.}
A kernel $\kernel$  is called strictly positive definite on $\paradom\subseteq\R^{\pdim}$ if for all $N\in\mathbb{N}$, all pairwise
distinct $Y_{N}=\left\{\vec{y}_1,\ldots,\vec{y}_N\right\}\subseteq\paradom$, and all
$\alpha\in\R^N\setminus\{0\}$ we have 
\begin{linenomath*}
\reviewerone{\begin{align}
\sum_{j=1}^N\sum_{k=1}^N \alpha_j
\alpha_k\,\kernel (\vec{y}_j,\vec{y}_k) > 0.
\end{align}}
\end{linenomath*}
Note that it is possible to construct a RKHS $\mathcal{H}_{\kernel}$ from a given strictly positive definite kernel function $\kernel$.
We then call $\mathcal{H}_{\kernel}$ the native space of $\kernel$, see \cite{Wendland2004}.

\subsection{Best approximation and regression in RKHS}\vspace*{-0.5em}
\label{sec:bestApproximation}
Now, let a strictly positive definite kernel function $\kernel$ with its
associated native space $\mathcal{H}_{\kernel}(\paradom)$ be given.
For the finite set of collocation points $Y_{\Nsamp}$ as in \eqref{collocpoints} we evaluate the kernel function $\kernel$ in these points. Thereby we introduce a finite-dimensional subspace
\begin{linenomath*}
\reviewerone{\begin{equation}
\mathcal{H}_{Y_{\Nsamp}} := \spn \{\kernel(\cdot,\vec{y}_i)|\vec{y}_i\in Y_{\Nsamp}\}\subset \mathcal{H}=\mathcal{H}_{\kernel}
\end{equation}}
\end{linenomath*}
of the reproducing kernel Hilbert space.  
Our aim is to approximate a given function $f\in \mathcal{H}_{\kernel}$ by a function in $\mathcal{H}_{Y_{\Nsamp}}$, i.e.
\begin{linenomath*}
\reviewerone{\begin{equation}
f(\vec{y}) \approx \sum_{i=1}^{\Nsamp} \alpha_i \kernel(\vec{y},\vec{y}_i)\,.
\end{equation}}
\end{linenomath*}
It is well-known that, if we only consider evaluations of $f$ in the collocation points $Y_{\Nsamp}$, we get the best approximation of $f$ with respect to the native space norm by computing coefficients $\vec{\alpha}:=(\alpha_1,\ldots, \alpha_{\Nsamp})^\top$ with
\begin{linenomath*}
\reviewerone{\begin{equation}
A_{\kernel, Y_{\Nsamp}} \vec{\alpha} = \vec{f}\,,
\end{equation}}
\end{linenomath*}
where $\vec{f}$ is the data vector $(f(\vec{y}_1),\ldots , f(\vec{y}_{\Nsamp}))^\top$. Here, the matrix $A_{\kernel,Y_{\Nsamp}}$ is given by
\begin{linenomath*}
\begin{equation}\label{eq:interpolationMatrix}
A_{\kernel,Y_{\Nsamp}} := 
\left(\begin{matrix}
\kernel(\vec{y}_1,\vec{y}_1) & \hdots & \kernel(\vec{y}_1,\vec{y}_{\Nsamp})\\
\vdots & \ddots & \vdots\\ 
\kernel(\vec{y}_{\Nsamp},\vec{y}_1) & \hdots & \kernel(\vec{y}_{\Nsamp},\vec{y}_{\Nsamp}) 
\end{matrix}\right)\,,
\end{equation}
\end{linenomath*}
where the strict positive definiteness of $\kernel$ ensures invertibility of this system.
We can build the Lagrange basis of $\mathcal{H}_{Y_{\Nsamp}}$ by setting
\begin{linenomath*}
\begin{align}\label{Lbasis}
 \left(L_{1}(\vec{z}),\dots,L_{\Nsamp}(\vec{z}) \right)^{\top}:=A^{-1}_{\kernel,Y_{\Nsamp}} \left( \kernel(\vec{z},\vec{y}_{1}), \dots,  \kernel(\vec{z},\vec{y}_{\Nsamp})\right)^{\top}. 
\end{align}
\end{linenomath*}
This leads to an approximation of $f$ in terms of the Lagrange basis of the form
\begin{linenomath*}
\begin{align}\label{eq:StochCollocInterpolationCondition}
f(\vec{y}) \approx \sum_{j=1}^{\Nsamp} L_{j}(\vec{y}) f(\vec{y}_{j}).
\end{align}
\end{linenomath*}
Once the Lagrange basis is determined, this is a favorable way to represent interpolation.

From a practical point of view, the matrix $A_{\kernel,Y_{\Nsamp}}$ may in general become 
ill-conditioned, depending on the kernel and the choice and number of collocation points.
One approach to tackle this issue is the introduction of a regularization. 
A standard technique is the Tikhonov
regularization \reviewerone{\cite{Tikhonov1977}}. It involves replacing the original linear system by
\begin{linenomath*}\vspace*{-1.5em}
\begin{align}\label{eq:regularizedKernelInterpolationLinearSystem}
(A_{\kernel,Y_{\Nsamp}}+\epsilon_{reg}\,\, I_{\Nsamp})  \vec{\alpha}&= \vec{f}\,,
\end{align}
\end{linenomath*}
with $I_N$ the identity matrix. 
This regularization reduces the condition number of $A_{\kernel,Y_{\Nsamp}}$, but
introduces a new error of the order of the regularization parameter
$\epsilon_{reg}$. Moreover, this regularization also accounts for the fact that we never deal with the true quantities of interest, e.g. \reviewerone{$f=u$ but always with some numerical approximations of them}. Hence, the regularized regression approach is favorable compared to classical interpolation. Nevertheless, we assume the data to be almost precise and hence $\epsilon_{reg}$ will be chosen very small, for example $\epsilon_{reg}=10^{-12}$ as in Section \ref{sec:NumResBenchmarkSetup}.

\reviewerone{An alternative to Tikhonov regularization is regularization by a truncated singular value decomposition (TSVD) \cite{tsvd}. In this approach, a singular value decomposition of matrix $A_{\kernel,Y_N}$ is computed. The SVD is truncated for singular values below a given magnitude. Then, instead of solving the linear system  $A_{\kernel, Y_{\Nsamp}} \vec{\alpha} = \vec{f}$, the pseudo-inverse of the truncated matrix $A_{\kernel,Y_N}$ is applied to $\vec{f}$. In Section~\ref{sec:ModelsProblemsWithAnalyticSolution}, we briefly compare numerical results with Tikhonov and TSVD regularization for a model problem.}

\subsection{Estimation of stochastic moments}\label{sec:stochasticMomentEstimation}
\reviewerone{The approximation of the stochastic moment $\E{u}$ needs the evaluation of an integral. To this end, conventional numerical quadrature methods like MC, QMC or sparse grids, etc. could be employed. Here, we restrict ourselves to kernel-based interpolatory quadrature rules.
We have} 
\begin{linenomath*}
\reviewerone{\begin{align}
\E{u}(\vec{x},t)=\int_{\paradom} u(\vec{y},\vec{x},t)\rho(\vec{y})\,d\vec{y} \approx \sum_{i=1}^{\Nsamp} u(\vec{y}_i,\vec{x},t) \int_\paradom L_{i}(\vec{y}) \rho(\vec{y}) d\vec{y},
\end{align}}
\end{linenomath*}
with the Lagrange basis from \eqref{Lbasis}. Moreover, we can reduce the integral
\begin{linenomath*}
\begin{align}
\int_\paradom L_i(\vec{y})\rho(\vec{y}) d\vec{y}&=\sum_{j=1}^{\Nsamp} c^{h}_{i,j} \reviewerone{\int_\paradom \kernel(\vec{y},\vec{y}_{j})\rho(\vec{y}) d\vec{y} \label{int1}}
\end{align}
\end{linenomath*}
to linear combinations of integrals over the kernel where the coefficients are given as $A^{-1}_{\kernel,Y_{\Nsamp}}=(c^{h}_{i,j})_{1\le i,j \le \Nsamp}$.
In order to approximate the kernel integral in \eqref{int1},
we employ here for the reason of simplicity full tensor-product Clenshaw-Curtis quadrature but also any other suitable method will do. Hence, the quadrature points are in general different from the sampling points. Note here that the quadrature of the kernel functions can be made arbitrarily precise and once the collocation points are designed, the integral in \eqref{int1} can be pre-computed anyway. Moreover, we do not consider the numerical costs for the quadrature.
This is due to the fact that a sampling point corresponds to a full three dimensional two--phase fluid simulation and hence is much more costly than a quadrature point.

\subsection{Sampling in $\paradom$ and choice of kernels}
\reviewerone{As collocation points $Y_{\Nsamp}$, we use low discrepancy points sets. They have the advantage that they allow to produce almost quasi-uniform points sets in higher dimensions. Since kernel based methods can work on arbitrarily scattered point sets, the specific choice is not crucial for the remainder of this article.} 
We use so-called \textit{quasi-random} sequences $Z_{QMC}:=\left\{\vec{z}_i\right\}_{i=1}^{\Nsamp}$ to sample collocation points. Quasi-random sequences
are designed to have a small discrepancy \cite{Caflisch1998}.
Here, to be precise, we employ multi-dimensional Halton sequences. 
The $i$th element of a
$\pdim$-dimensional Halton sequence in $[0,1)^{\pdim}$ is given
as
\begin{linenomath*}
\reviewerone{\begin{equation}
\vec{z}_{i} = (\psi_{R_1}(i),\psi_{R_2}(i),\ldots,
\psi_{R_{\pdim}}(i)),\end{equation}}
\reviewerone{\begin{equation}\psi_{R}(i) = \sum_{j=1}^\infty a_j(i) R^{-j},\,\,\,\,\,\, \mbox{with}\
\,\,\,\,i= \sum_{j=1}^{\infty} a_j(i) R^{j-1}.
\end{equation}}
\end{linenomath*}
Here, the sums are all finite and the $R_{k}$ have to be coprime integers, and the $a_{{j}}(i)$ just denote the $R$-adic representation of $i$. \reviewerone{We employ the Halton points here, since we work in moderately high dimensions. We are well aware of the fact that Halton points might deteriorate in high dimensions and scrambling might be needed then, see \cite{Chietal,QMCbook,WangSloan}. Moreover, we stress that we work with a kernel based and hence mesh free method. This implies that our methods works on all scattered point sets and only the theoretical convergence analysis is affected by the distribution properties of the point set.}

Our approach was so far described for any suitable RKHS. The final question is now what types of RKHS we want to invoke and what the associated kernel $\kernel:\paradom \times \paradom \to \R$ will be. To this end, we restrict ourselves to
radial basis functions/radial kernels. 
They are given by 
\begin{linenomath*}
\reviewerone{\begin{equation}
\kernel (\vec{y},\vec{y}^\prime) := 
\varphi(||\vec{y}-\vec{y}^\prime||_2)\,,
\end{equation}}
\end{linenomath*}
with an appropriate function $\varphi:\R^{\geq 0}\rightarrow \R$. \begin{table}
\begin{center}
\begin{tabular}{c|c|c}\hline
kernel  & definition & Sobolev space\\\hline\hline
Wendland & $\kernel_{\pdim,k} (\vec{y},\vec{y}^\prime) :=
\varphi_{\pdim,k}(\|\vec{y}-\vec{y}^\prime\|)\,,$   & $H^{\pdim/2+k+1/2}(\R^{\pdim})$\\\hline
Mat\'ern & $\kernel_{\beta}(\vec{y},\vec{y}^\prime) :=
\frac{K_{\beta-\frac{\pdim}{2}}(\|\vec{y}-\vec{y}^\prime\|)\|\vec{y}-\vec{y}^\prime\|^{\beta-\frac{\pdim}{2}}}{2^{\beta
-1}\mathsf{\paradom}(\beta)}$ & $ W^{\beta,2}(\R^{\pdim})$\\\hline
Gaussian & $k_\epsilon (\vec{y},\vec{y}^\prime) :=
\varphi_\epsilon(\|\vec{y}-\vec{y}^\prime\|) := e^{-\epsilon^2
\|\vec{y}-\vec{y}^\prime\|^2}$ & n/a\\
\end{tabular}
\end{center}
\caption{\label{tab:kernels}The above radial kernel functions are used in this work to build native (approximation) spaces. The resulting native spaces are well-known Sobolev spaces.}
\end{table}

\noindent
There are several examples of such radial kernels. Table~\ref{tab:kernels} summarizes the kernel functions used in this work and the Sobolev spaces that are equal to their native space. 

Wendland kernels $\kernel_{d,k}$ \cite{Wendland2004} are compactly supported functions with some minimality properties for the degrees of the polynomials involved in their construction. They are positive definite and it holds $\varphi_{\pdim,k}\in C^{2k}(\R)$. For their associated $\kernel_{\pdim,k}$ for $k=0,1$ we e.g.~have 
\begin{linenomath*}
\reviewerone{\begin{equation}
\varphi_{\pdim,0}(r) = (1-r)_+^{\lfloor \pdim/2\rfloor +1}, \quad
\varphi_{\pdim,1}(r) = (1-r)_+^{\ell+1}[(\ell +1)r+1]\,,
\end{equation}}
\end{linenomath*}
with $\ell:=\lfloor \pdim/2\rfloor+k+1$ and the notation $(r)_+ = \left\{\begin{array}{ll}
	r & \mbox{if}\ r\geq 0\,,\\
	0 & \mbox{if}\ r< 0\,.	
\end{array}\right.$ 

The Mat\'ern kernels $\kernel_{\beta}$ with $\beta>\frac{d}{2}$ use in their definition $K_\nu$ and $\mathsf{\paradom}$, which are the modified Bessel function of the second kind of order $\nu$ and the Gamma function, see also \cite{Wendland2004}. These kernels are strictly positive definite as long as $d < 2\beta$. According to \cite[Section~4.4]{Fasshauer2007}, we have for special choices of
parameter $\beta$ simplified representations of the Mat\'ern kernel
function (up to a dimension-dependent scaling constant), e.g.~
\begin{linenomath*}
\reviewerone{\begin{equation}
\kernel_{\frac{\pdim+1}{2}}(\vec{y},\vec{y}^\prime) :=
e^{-\|\vec{y}-\vec{y}^\prime\|}, \quad\kernel_{\frac{\pdim+3}{2}}(\vec{y},\vec{y}^\prime) :=
\left(1+\|\vec{y}-\vec{y}^\prime\|\right) e^{-\|\vec{y}-\vec{y}^\prime\|}\,.
\end{equation}}
\end{linenomath*}
The last example is the well-known Gaussian kernel with $\epsilon\in\R^+$ as scaling parameter. It is special in the sense that its native space is contained in every Sobolev space.

\section{Numerical results}
\label{sec:numericalResults}
In the following, we will discuss our method for \reviewerone{two analytic test cases and} several two-phase flow problems. After the introduction of the setup, empirical convergence results will be presented for kernel-based stochastic collocation using different kernel functions. At the end of this section, we also briefly compare results from our kernel-based approach with a sparse grid approach for stochastic collocation.

\subsection{Setup}\vspace*{-1em}
\label{sec:NumResBenchmarkSetup}
In this section, \reviewerone{we first consider the approximation of means of two analytic test functions, cf.~Section~\ref{sec:ModelsProblemsWithAnalyticSolution}.
Thereafter,} we compute approximations to the field
$\E{u^{h}_1}(\vec{x},T)$ and the scalar $\E{c^{h}_2}(T)$ for fixed $T$.
The deterministic evaluations of the quantities of interest $u^{h}_1$ and $c^{h}_2$ are
done using NaSt3DGPF for each $\vec{y}_{i}\in \paradom$. Since this solver uses a finite volume/finite difference method on
a staggered grid, we obtain scalar-valued, space-dependent fields as grid functions
on a regular grid with $\Ngrid^{h}:=\Ngrid$ points. Here, the parameter $h$, which will be dismissed in the remainder, indicates that the $\Ngrid$ points correspond to a  grid in space with uniform meshsize $h$ and an accordingly chosen CFL conforming time-step.  
In case of the quantity of interest $u^{h}_1$, we thus actually
compute point-wise for every grid point the first moment $\E{u^{h}_1}(T)\in\R^{\Ngrid}$.
The quantity $\E{c^{h}_2}(T)$ is a single number.

To approximate these quantities with respect to the stochastic space, we use
the RBF kernel-based stochastic
collocation method. If not stated otherwise, the applied kernel functions are
the Gaussian kernel $\kernel_\epsilon$ with scaling parameter $\epsilon=1.0$, compactly supported Wendland kernels $\kernel_{\pdim,k}$ with smoothness parameters
$k=0,1,2,3$ and appropriate dimensionality $\pdim$ and the Mat\'ern
kernel $\kernel_\beta$ with parameter $\beta=\frac{\pdim+3}{2}$. Remember that
we call the dimension of the stochastic space $\pdim$, thus $\paradom\subset\R^{\pdim}$.
The regularization parameter,
cf.~Section~\ref{sec:bestApproximation}, is set to $\epsilon_{reg}=10^{-12}$.
The kernel interpolation problem is solved by direct LU factorization. 
Radial basis functions are isotropic, by
standard, thus the norm involved in their construction is the (scaled)
Euclidean distance $\|\cdot\|:=\kernelscale\|\cdot\|_2$, with a default of
$\kernelscale=1.0$.
\reviewerone{In Section~\ref{sec:ModelsProblemsWithAnalyticSolution}, we briefly compare different
choices of $\kernelscale$ for a model problem.}
The $\Nsamp$ collocation points are generated from a Halton sequence of appropriate
dimension. 

Quadrature is carried out by a full tensor product rule constructed by univariate Clenshaw-Curtis quadrature rules with $2^{l_q-1}+1$ nodes, if not indicated differently. The quadrature-level is usually $l_q=7$.
Sparse grid quadrature rules are used whenever the
dimensionality of the stochastic space would lead to prohibitive
computational run-times and memory requirements.

In the empirical convergence studies, the reference solutions will be computed on an extremely fine grid.
We will use a subscript $H$ to denote the reference solutions.
Moreover, in order to visualize the numerical error, we approximatively compute the expectation of the reference solution as 
\begin{linenomath*}
\begin{eqnarray}\label{reffirstmom}
\E{u_{i}^{H}}(\vec{x},T):=\int_{\paradom} u^{H}_{i}(\vec{y},\vec{x},T)\rho(\vec{y})d\vec{y}\approx \sum_{n=1}^{N_{\max}}u_{i}^{H}(\vec{\tilde{y}}_{n},\cdot,T) \int_{\paradom} \tilde{L}_{n}(\vec{y})\rho(\vec{y})d\vec{y},
\end{eqnarray}
\end{linenomath*}
where $\tilde{Y}_{N_{\max}} \subset \paradom$ denotes a fine sampling set, i.e,. $N_{\max}>> N$. Here, $\tilde{L}$ is computed based on a possibly different kernel compared to the approximation.

\subsection{\reviewerone{Problems with analytic solution}}
\label{sec:ModelsProblemsWithAnalyticSolution}

Before we discuss numerical results for the challenging two-phase Navier-Stokes
application problem, we briefly analyse the properties of the kernel-based stochastic collocation
for two representative test cases with analytic solution. \crcom{To this end, we restrict ourselves to a simple model problem.}

The first test case is an elliptic PDE with random coefficient similar to \cite{Tuminaro2011}.
The parameter space is $(\paradom, \mathcal{B},\rho\,d\vec{y})$ with $\paradom\subset\R$ and stems
from a one-term \KL expansion of a random coefficient. We want to approximate $\E{u}$ with 
$u:\paradom\times [-0.5,0.5]^2\rightarrow\R$ being the solution of
\begin{linenomath*}
\begin{eqnarray}\label{testcase1}
-\nabla\cdot(a(\vec{y},\vec{x})\nabla u(\vec{y},\vec{x})) =& f(\vec{x}) &
\mbox{in}\ \paradom\times (-0.5,0.5)^2\,,\\ \label{testcase2}
u(\vec{y},\vec{x}) =& 0  & \mbox{on}\ \paradom\times \partial(-0.5,0.5)^2\,, 
\end{eqnarray}
\end{linenomath*}
with the random diffusion coefficient
$a(\vec{y},\vec{x})=1 + \sigma \frac{1}{\pi^2} y_1
\cos\left(\frac{\pi}{2}\left(x_1^2+x_2^2\right)\right),$
thus we have $\vec{y}=y_1$. The right-hand side term is given as
\begin{linenomath*}
\begin{align*}
f(\vec{y},\vec{x}) =& 32 \left(1+\sigma + \frac{y_1
\cos(\frac{1}{2}\pi(x_1^2+x_2^2))}{\pi^2}\right)
e^{-y_1^2} \left(x_2^2 -\frac{1}{2} + x_1^2 \right)\\
& -\frac{32}{\pi} y
\sin\left(\frac{1}{2}\pi (x_1^2+x_2^2)\right) \left(x_1^2 e^{-y_1^2}\left(x_2^2-\frac{1}{4}\right) + x_2^2
e^{-y_1^2}\left(x_1^2-\frac{1}{4}\right) \right)\,.
\end{align*}
\end{linenomath*}
With this construction, it is possible to derive an
exact solution of the parametric PDE problem as
\begin{linenomath*}
\begin{equation}u(\vec{y},\vec{x})=16\,e^{-y_1^2}
\left(x_1^2-\frac{1}{4}\right)\left(x_2^2-\frac{1}{4}\right)\,.
\end{equation}
\end{linenomath*}
The variable $y_1$ corresponds to the random variable
$Y_1(\omega)\sim\mathcal{U}(-\sqrt{3},\sqrt{3})$, thus we employ the 
density function $\rho(\vec{y})=\frac{1}{2\sqrt{3}}$.
It is possible to derive the exact mean as
$\E{u}=\frac{1}{6}\erf\left(\sqrt{3}\right)\sqrt{3}\sqrt{\pi}\left(16 x_1^2
x_2^2 - 4 x_1^2 -4 x_2^2 + 1\right).$


\begin{figure}
\begin{center}\vspace*{-1em}
\hfill\scalebox{0.75}{\begin{tikzpicture}
\begin{loglogaxis}[
	height=7.8cm,
    width=7.8cm, 
	xlabel=\# collocation points,
	ylabel=absolute error,
	legend style={legend pos=north east,font=\tiny},
	ymax=1.0e+1
	]
\addplot
table[x=collocationpoints,y=gaussian]{error_mean_poisson1d_kernel_overview_new.csv};
\addplot
table[x=collocationpoints,y=wendland0]{error_mean_poisson1d_kernel_overview_new.csv};
\addplot
table[x=collocationpoints,y=wendland1]{error_mean_poisson1d_kernel_overview_new.csv};
\addplot
table[x=collocationpoints,y=wendland2]{error_mean_poisson1d_kernel_overview_new.csv};
\addplot
table[x=collocationpoints,y=wendland3]{error_mean_poisson1d_kernel_overview_new.csv};
\addplot
table[x=collocationpoints,y=matern]{error_mean_poisson1d_kernel_overview_new.csv};
\legend{
Gauss. $\kernel_\epsilon$,
Wend. $\kernel_{1,0}$,  
Wend. $\kernel_{1,1}$,
Wend. $\kernel_{1,2}$,
Wend. $\kernel_{1,3}$,
{Mat\'ern $\kernel_{\frac{1+3}{2}}$}}
\addplot[black] coordinates {
                (1000, 1.0e-7)
               (8000, 1.5625e-9)
        };
\draw (axis cs:8000,1.5625e-9) node[right]{{\small 2}};
\addplot[black] coordinates {
                (1000, 1.0e-9)
               (8000, 1.953125e-12)
        };
\draw (axis cs:8000,1.953125e-12) node[right]{{\small 3}};
\addplot[black] coordinates {
				(100, 1.0e-9)
				(800, 3.0518e-14)
		};
\draw (axis cs:800,3.0518e-14) node[right]{{\small 5}};
\end{loglogaxis}
\end{tikzpicture}}\hfill\scalebox{0.75}{\begin{tikzpicture}
\begin{loglogaxis}[
	height=7.8cm,
    width=7.8cm, 
	xlabel=\# collocation points,
	ylabel=absolute error,
	legend style={legend pos=south west,font=\tiny},
	ymax=1.0e+1
	]
\addplot
table[x=collocationpoints,y=m4]{error_mean_poisson1d_wendland3_scalings_new.csv};
\addplot
table[x=collocationpoints,y=m2]{error_mean_poisson1d_wendland3_scalings_new.csv};
\addplot
table[x=collocationpoints,y=m1]{error_mean_poisson1d_wendland3_scalings_new.csv};
\addplot
table[x=collocationpoints,y=p0]{error_mean_poisson1d_wendland3_scalings_new.csv};
\addplot
table[x=collocationpoints,y=p1]{error_mean_poisson1d_wendland3_scalings_new.csv};
\addplot
table[x=collocationpoints,y=p2]{error_mean_poisson1d_wendland3_scalings_new.csv};
\addplot
table[x=collocationpoints,y=p4]{error_mean_poisson1d_wendland3_scalings_new.csv};
\legend{
$\zeta=10^{-4}$,
$\zeta=10^{-2}$,
$\zeta=10^{-1}$,
$\zeta=10^{0}$,
$\zeta=10^{1}$,
$\zeta=10^{2}$,
$\zeta=10^{4}$
}
\end{loglogaxis}
\end{tikzpicture}}\hfill\,
\end{center}\vspace*{-1em}
\caption{\label{fig:ResultsErrorMeanAndSecondMomentPoisson1D}\reviewerone{Error
convergence analysis for the mean approximation in the random-coefficient Poisson problem test case \eqref{testcase1} \& \eqref{testcase2} with $D=1$.
Left: Comparison of different kernel functions. Right: Comparison of different scalings $\zeta$ for the fixed kernel function $\kernel_{1,3}$.}}
\end{figure}
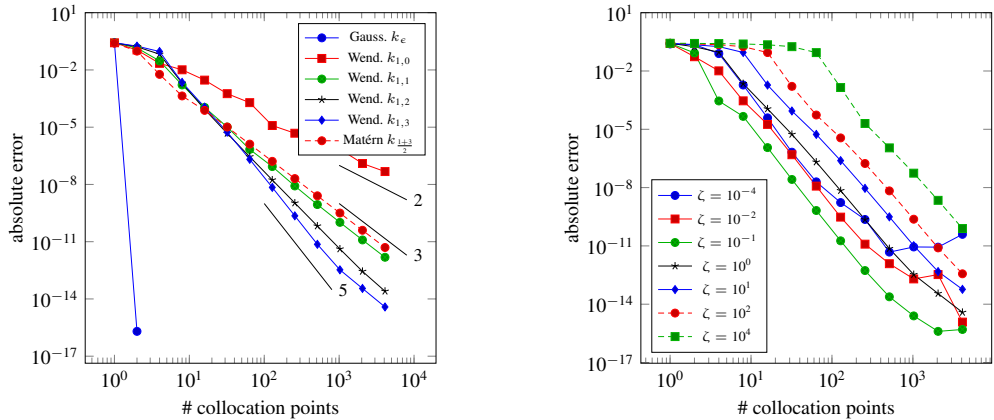

The convergence study for the mean is presented on the left-hand side of Figure~\ref{fig:ResultsErrorMeanAndSecondMomentPoisson1D}.
It shows absolute errors of the mean solution field with respect to the exact solution measured
in a discrete $l_2$ norm and compares different kernel choices.
The regularization parameter $\epsilon_{reg}$ is set to $10^{-15}$ for
the Gaussian kernel. All other parameters remain at the previously defined values, see Section \ref{sec:NumResBenchmarkSetup}.
Approximation with the Gaussian kernel leads \crcom{to an almost perfect solution, i.e., 
machine precision with only two collocation points, which is due to the specific choice of the 
unknown function $u$.} All other kernels give algebraic convergence rates with
measured approximate orders 2, 3, 4 and 5 for Wendland kernels with $k=0,1,2,3$
and third-order convergence for that specific choice of a Mat\'ern kernel.
Convergence saturates between machine accuracy and the size of the regularization parameter, which, in contrast to the Gaussians, is set to $\epsilon_{reg}=10^{-12}$ for
these kernels. 

On the right-hand side of Figure~\ref{fig:ResultsErrorMeanAndSecondMomentPoisson1D}, we
repeat the same convergence study. However, we now fix the Wendland kernel $\kernel_{1,3}$
and change the scaling $\kernelscale$. The results indicate that this scaling has
a substantial influence on the convergence of the method. In the given results, a choice of $\kernelscale$, which is too large, leads to a longer pre-asymptotic regime. A choice
of $\kernelscale$, which is too small, leads to convergence issues. For the rest of this paper, we will manually
optimize this parameter.

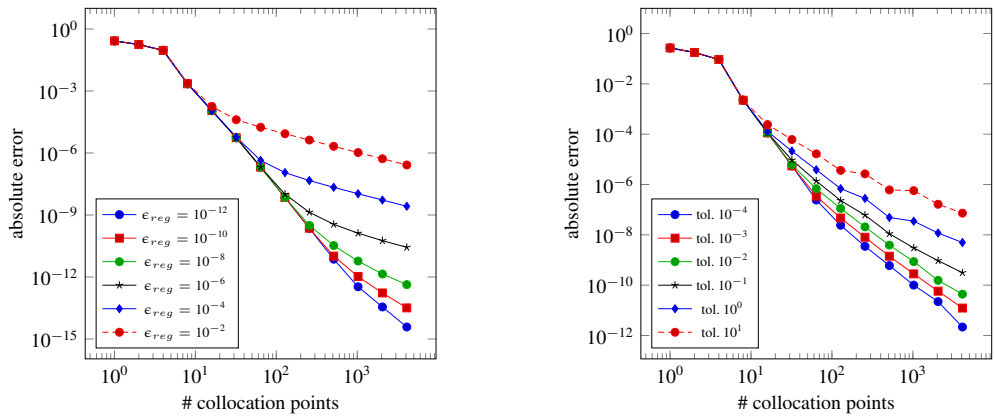
\begin{figure}
\begin{center}\vspace*{-1em}
\hfill\scalebox{0.75}{\begin{tikzpicture}
\begin{loglogaxis}[
	height=7.8cm,
    width=7.8cm, 
	xlabel=\# collocation points,
	ylabel=absolute error,
	legend style={legend pos=south west,font=\tiny},
	ymax=1.0e+1
	]
\addplot
table[x=collocationpoints,y=eps12]{error_mean_poisson1d_wendland3_regularizations_new.csv};
\addplot
table[x=collocationpoints,y=eps10]{error_mean_poisson1d_wendland3_regularizations_new.csv};
\addplot
table[x=collocationpoints,y=eps08]{error_mean_poisson1d_wendland3_regularizations_new.csv};
\addplot
table[x=collocationpoints,y=eps06]{error_mean_poisson1d_wendland3_regularizations_new.csv};
\addplot
table[x=collocationpoints,y=eps04]{error_mean_poisson1d_wendland3_regularizations_new.csv};
\addplot
table[x=collocationpoints,y=eps02]{error_mean_poisson1d_wendland3_regularizations_new.csv};
\legend{
$\epsilon_{reg}=10^{-12}$,
$\epsilon_{reg}=10^{-10}$,
$\epsilon_{reg}=10^{-8}$,
$\epsilon_{reg}=10^{-6}$,
$\epsilon_{reg}=10^{-4}$,
$\epsilon_{reg}=10^{-2}$,
}
\end{loglogaxis}
\end{tikzpicture}}\hfill\scalebox{0.75}{\begin{tikzpicture}
\begin{loglogaxis}[
	height=7.8cm,
    width=7.8cm, 
	xlabel=\# collocation points,
	ylabel=absolute error,
	legend style={legend pos=south west,font=\tiny},
	ymax=1.0e+1
	]
\addplot
table[x=collocationpoints,y=m4]{error_mean_poisson1d_wendland3_svd_regularizations_new.csv};
\addplot
table[x=collocationpoints,y=m3]{error_mean_poisson1d_wendland3_svd_regularizations_new.csv};
\addplot
table[x=collocationpoints,y=m2]{error_mean_poisson1d_wendland3_svd_regularizations_new.csv};
\addplot
table[x=collocationpoints,y=m1]{error_mean_poisson1d_wendland3_svd_regularizations_new.csv};
\addplot
table[x=collocationpoints,y=p0]{error_mean_poisson1d_wendland3_svd_regularizations_new.csv};
\addplot
table[x=collocationpoints,y=p1]{error_mean_poisson1d_wendland3_svd_regularizations_new.csv};
\legend{
tol. $10^{-4}$,
tol. $10^{-3}$,
tol. $10^{-2}$,
tol. $10^{-1}$,
tol. $10^{0}$,
tol. $10^{1}$,
}
\end{loglogaxis}
\end{tikzpicture}}\hfill\,
\end{center}\vspace*{-1em}
\caption{\label{fig:RegularizationComparisonPoisson1D}\reviewerone{Comparison of the Tikhonov regularization (left) and the TSVD regularization
(right) for the random-coefficient Poisson problem test case. The identical Wendland kernel $\kernel_{1,3}$ is applied.}}
\end{figure}

Figure~\ref{fig:RegularizationComparisonPoisson1D} shows results for the use of Tikhonov regularization
on the left-hand side and for the use of the truncated SVD on the right-hand side. As before, we approximate
the mean, fixing the kernel function $\kernel_{1,3}$. However, we do variations in either the regularization
parameter $\epsilon_{reg}$ (for Tikhonov) or in the dropping tolerance for singular values (for TSVD). Both
approaches show similar regularization properties, i.e.~stronger regularization leads to a higher total error
in the approximation. Throughout the rest of this work, we will use Tikhonov regularization, since this
approach seemed to be numerically more stable for smaller regularization. 

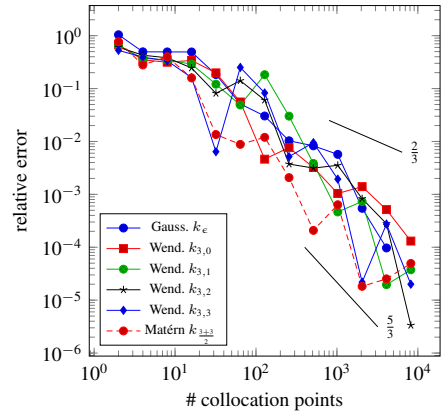
\begin{figure}
\begin{center}\scalebox{0.75}{\begin{tikzpicture}
\begin{loglogaxis}[
	height=7.8cm,
    width= 7.8cm,
	xlabel=\# collocation points,
	ylabel=relative error,
	legend style={legend pos=south west,font=\tiny}
	]
\addplot
table[x=collocationpoints,y=gaussian3d]{error_mean_g_function_kernel_overview.csv};
\addplot
table[x=collocationpoints,y=wendland3d0]{error_mean_g_function_kernel_overview.csv};
\addplot
table[x=collocationpoints,y=wendland3d1]{error_mean_g_function_kernel_overview.csv};
\addplot
table[x=collocationpoints,y=wendland3d2]{error_mean_g_function_kernel_overview.csv};
\addplot
table[x=collocationpoints,y=wendland3d3]{error_mean_g_function_kernel_overview.csv};
\addplot
table[x=collocationpoints,y=matern]{error_mean_g_function_kernel_overview.csv};
\legend{Gauss. $\kernel_\epsilon$, Wend.
$\kernel_{3,0}$, Wend.  $\kernel_{3,1}$, Wend.
$\kernel_{3,2}$, Wend.  $\kernel_{3,3}$, 
{Mat\'ern $\kernel_{\frac{3+3}{2}}$}}
\addplot[black] coordinates {
                (800, 2.5e-2)
               (6400, 6.25e-3)
        };
\draw (axis cs:6400,6.25e-3) node[right]{{\small $\frac{2}{3}$}};
\addplot[black] coordinates {
                (400, 1.0e-4)
               (3200, 3.125e-6) 
        };
\draw (axis cs:3200,3.125e-6) node[right]{{\small $\frac{5}{3}$}};
\end{loglogaxis}
\end{tikzpicture}}
\end{center}\vspace*{-1em}
\caption{\label{fig:ResultsErrorMeanGfunction}\reviewerone{Error
convergence analysis for the mean approximation comparing different kernel functions
for the function \eqref{davis}.}}
\end{figure}

As a higher-dimensional test problem with limited smoothness, we further consider
the function
\begin{equation}\label{davis} u(\vec{y}) =
\prod_{m=1}^{D} \frac{|4 y_m -2| + a_m}{1 + a_m}\,,\end{equation} cf.~\cite{Davis1984}.
Here, we approximate the exact mean $\E{u} = 1$ with the choices $a_m = \frac{m-2}{2}$, $\rho(\vec{y})=1$. Moreover, the $y_m$ are realizations of independent random variables $Y_m(\omega)\sim\mathcal{U}(0,1)$.

  
In the numerical results, the scaling of the
Gaussian kernel is set to $\epsilon=2.0$.
Furthermore the approximation by the Gaussian kernel is regularized with a
regularization parameter $\epsilon_{reg}=10^{-8}$. We set $D=3$. 
Figure~\ref{fig:ResultsErrorMeanGfunction}
shows the error behavior of the mean of this model problem
for different kernels. We observe convergence rates which are better than that for Monte Carlo and even
better than that for quasi Monte Carlo methods. However, all kernel functions show qualitatively
identical results. This suggests that all results are affected by the limited smoothness
of the function \eqref{davis}, in the first place.

\subsection{Bubble flow under random volume force}\vspace*{-1em}
\label{sec:ModelsKLbasedRandomVolumeForceInBubbleFlow}
Our first application for the two-phase Navier--Stokes equations is a rising air bubble in water under a random volume force field with known covariance spectrum. 
The random force field is approximated by a
truncated \KL expansion. All other parameters are deterministic.
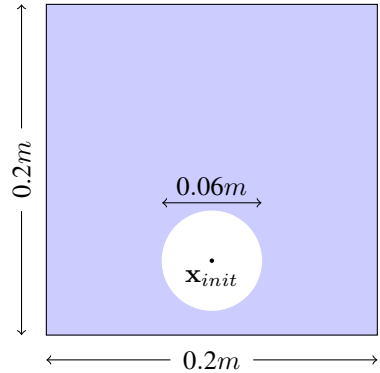
\begin{figure}[t]
\begin{center}\vspace*{0.5em}
\scalebox{1.0}{
\begin{tikzpicture}[scale=1.1]
\filldraw[lightblue] (0,0) rectangle (4,4);
\draw (0,0) rectangle (4,4);
\filldraw[white] (2,0.9) circle (0.6);
\filldraw (2,0.9) circle (0.02) node[below]{$\vec{x}_{init}$};

\draw[<->] (1.4,1.6) -- (2.6,1.6);
\draw (2,1.8) node{$0.06m$};

\draw[<-] (0,-0.3) -- (1.5,-0.3);
\draw[->] (2.5,-0.3) -- (4,-0.3);
\draw (2,-0.3) node{$0.2m$};

\draw[<-] (-0.3,0) -- (-0.3,1.5);
\draw[->] (-0.3,2.5) -- (-0.3,4);
\draw (-0.3,2) node[rotate=90]{$0.2m$};
\end{tikzpicture}}
\end{center}\vspace*{-1em}
\caption{\label{fig:ModelRisingBubbleSetup}Two-dimensional side view of the
three-dimensional rising bubble application setup.}
\end{figure}
Figure~\ref{fig:ModelRisingBubbleSetup} outlines the basic setup of the
two-phase flow problem\footnote{Note that some quantities have physical units. We do not write them every time in order to keep the notation simple. We will assume the densities $\rho$ in \reviewertwo{$[kg \ m^{{-3}}]$}, the viscosities in $[kg \ m^{-1} s^{{-1}}]$, the surface tension in $[N m^{{-1}}]$ and the gravitational force $\vec{g}$ in $[m s^{{-2}}]$.}. The domain is given as $\bar{\dom}=[0,0.2]^3$ and the
fluid flow is computed up to $T=0.35$ seconds.
Since an air-water system shall be analyzed, we have the densities
$\rho_1=1000$ and $\rho_2=1$ and viscosities $\mu_1=1.002\cdot 10^{-3}$ and
$\mu_2=1.72\cdot 10^{-5}$.\footnote{The dimensionless Reynolds number is given as $Re=\frac{\rho v L}{\mu}$, where $v$ and $L$ are the characteristic velocity and the characteristic length. Here, the Reynolds number is in the range $Re\in [10^{4},10^{6}]$. The validity of the discretization in the regime presented here was shown in \cite{Croce2009}.}
The initial conditions are set to $\vec{u}^{h}_{0_i}=\left(0,0,0\right)^\top$. Moreover, the surface tension
coefficient is set to $\sigma_{}=0.0728$ which reflects the
parameter of a water-air interface. Boundary conditions of the velocity
field are described in terms of a boundary operator 
\begin{linenomath*}
\reviewerone{\begin{equation}
\mathcal{B}\vec{u}^{h}:=\left(\vec{u}^{h}\cdot\vec{n}^{h}_{\bnd}, \frac{\partial(\vec{u}^{h}\cdot\vec{s}^{h}_{\bnd})}{\partial\vec{n}^{h}_{\bnd}}, \frac{\partial(\vec{u}^{h}\cdot\vec{t}^{h}_{\bnd})}{\partial\vec{n}^{h}_{\bnd}}\right)^{\top}=(0,0,0)^{\top}=b_{\bnd}
\end{equation}}
\end{linenomath*}
with the notation as in \eqref{eq:ModelDeterministicRandomNSbc1}. Thus, infinite slip is assumed on the boundary. 
Furthermore, the initial position of the center of the bubble is $\vec{x}^{init}=(0.1, 0.06, 0.1)^\top$. The phase-wise
sub-domains $\dom^{h}_1,\dom^{h}_2$ are given by
\begin{linenomath*}
\reviewerone{\begin{equation}
\dom^{h}_i(t) = \Phi\left[(\dom^{h}_i)^0\right](t)\,.
\end{equation}}
\end{linenomath*}
Here, $\Phi$ describes the
transformation of the initial domains $(\dom^{h}_i)^0:=\dom^{h}_i(t=0)$ under fluid flow.
These initial domains are defined such that the gas phase is a sphere of radius $0.03\,m$, thus
\begin{linenomath*}
\reviewerone{\begin{equation}
(\dom^{h}_1)^0 := \left\{ \vec{x}\in\dom
\middle| \|\vec{x}-\vec{x}_{init}\| > 0.03\right\}\,,
\end{equation}}
\end{linenomath*}
\begin{linenomath*}
\reviewerone{\begin{equation}
(\dom^{h}_2)^0 := \left\{
\vec{x}\in\dom \middle| \|\vec{x}-\vec{x}_{init}\| < 0.03\right\}\,.
\end{equation}}
\end{linenomath*}
The initial free surface is $\bnd_f(t=0) =
\dom\setminus(\dom_1^0\cup \dom_2^0)$.
In order to model a stochastic/parametric volume force $g$, we consider a lognormal random \reviewerone{field}
$g_{L_{c}}$ with $\E{g_{L_{c}}}\equiv -9.81$ and
\begin{linenomath*}
\reviewerone{\begin{equation}
\cov{\log (g_{L_{c}} - (-9.81))}(x,x^\prime) =
e^{-\frac{\left(x-x^\prime\right)^2}{L_c^2}},
\end{equation}}
\end{linenomath*}
where the correlation length is assumed to be $L_c=2.0$.
The random \reviewerone{field} $g_{L_{c}}$ is approximated by a
truncated \KL expansion resulting into $g_{\pdim,L_{c}}$. Consequently,
\begin{linenomath*}
\begin{eqnarray}\label{eq:ModelBubbleKLexpansion}
\log\left(g_{L_c}(\vec{y},\vec{x})+9.81\right)\approx \log\left(g_{\pdim,L_c}(\vec{y},\vec{x})+9.81\right) := 1+ y_1
\left(\frac{\sqrt{\pi} L_c}{2}\right)^{1/2} + \sum_{m=2}^{\pdim}
\lambda_m\phi_m(x_2)y_m\,,
\end{eqnarray}
\end{linenomath*}
with truncation after $\pdim$ expansion terms. 
Here, the
$\left\{y_m\right\}_{m=1}^{\pdim}$ correspond to the independent random variables\linebreak
$\left\{Y_m(\vec{\omega})\right\}_{m=1}^{\pdim}$ with each $Y_m\sim
\mathcal{U}(-\sqrt{3},\sqrt{3})$, thus
$\rho(\vec{y})=\left(\frac{1}{2\sqrt{3}}\right)^{\pdim}$. Furthermore, the
eigenvalues and eigenfunctions are given as\vspace*{-1em}
\begin{linenomath*}
\reviewerone{\begin{equation}
\lambda_m :=
(\sqrt{\pi}L_c)^{1/2} \exp\left(\frac{-(\lfloor\frac{m}{2}\rfloor \pi L_c)^2}{8}\right),\,\, \phi_m(x_2):=\left\{\begin{array}{ll}\sin (\lfloor\frac{m}{2}\rfloor\pi x_2) &
\mbox{if}\ m\ \mbox{even},\\\cos(\lfloor\frac{m}{2}\rfloor\pi x_2) & \mbox{if}\
m\ \mbox{odd},\end{array}\right.
\end{equation}}
\end{linenomath*}
with $m>1$.
Finally, the volume force is set to\vspace*{-1em}
\begin{linenomath*}
\reviewerone{\begin{equation}
\vec{g}(\vec{y},\vec{x}):= \left(0, g_{\pdim,L_c}(\vec{y},x_2),
0\right)^\top\,.\vspace*{-1em}
\end{equation}}
\end{linenomath*}\begin{figure}
\begin{center}
\scalebox{1.06}{\scalebox{0.055}{\includegraphics{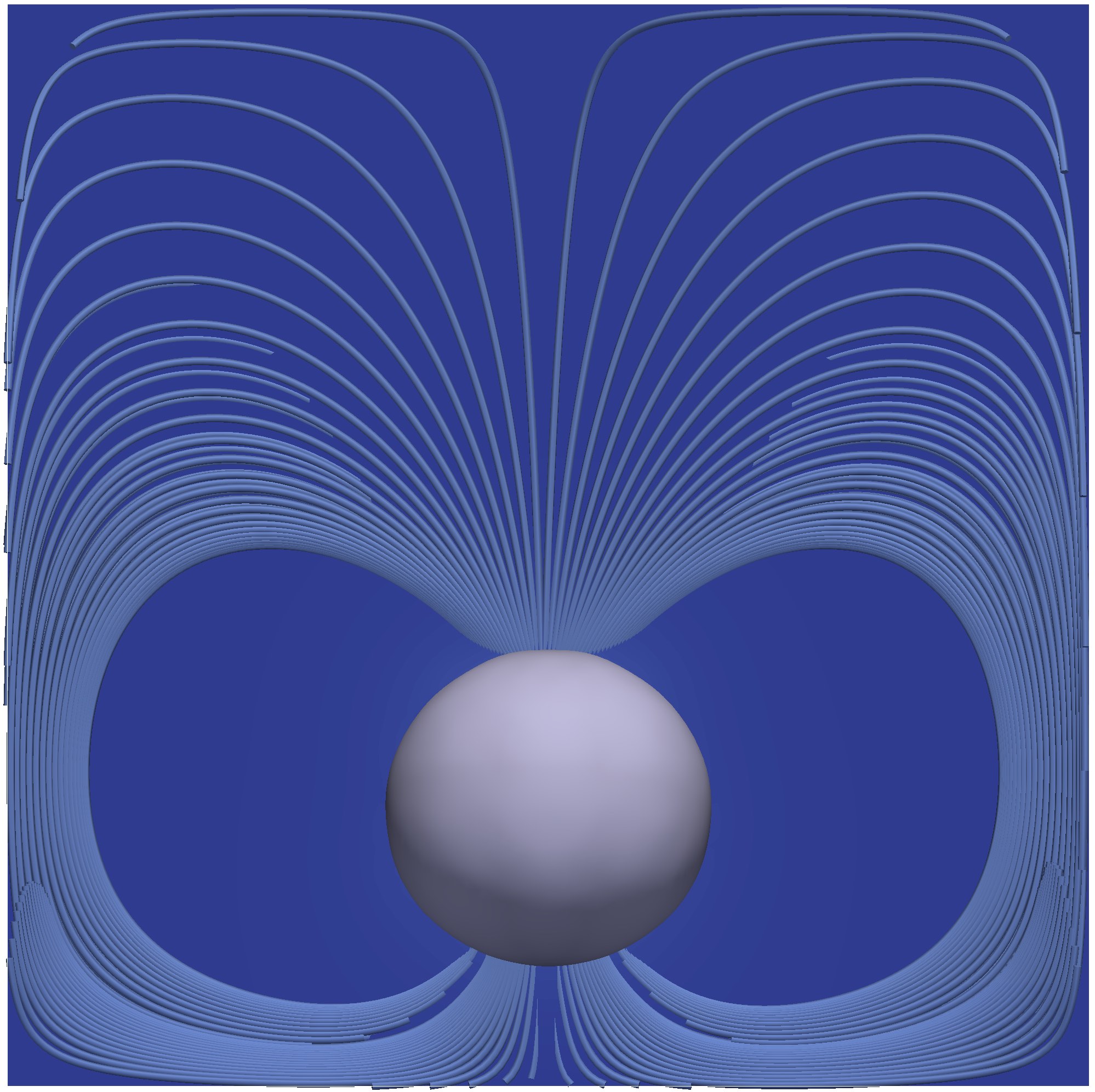}}\scalebox{0.055}{\includegraphics{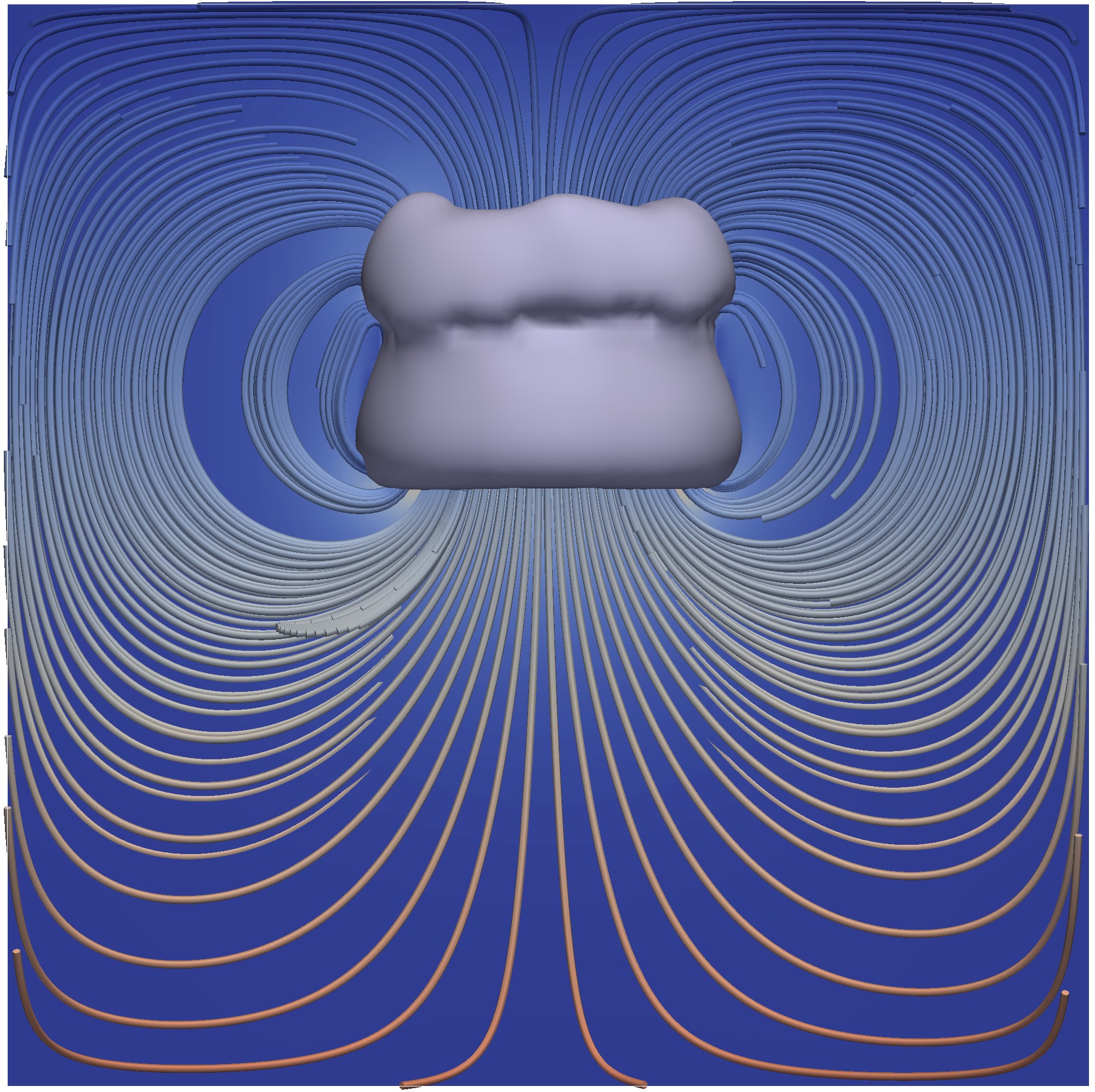}}\scalebox{0.055}{\includegraphics{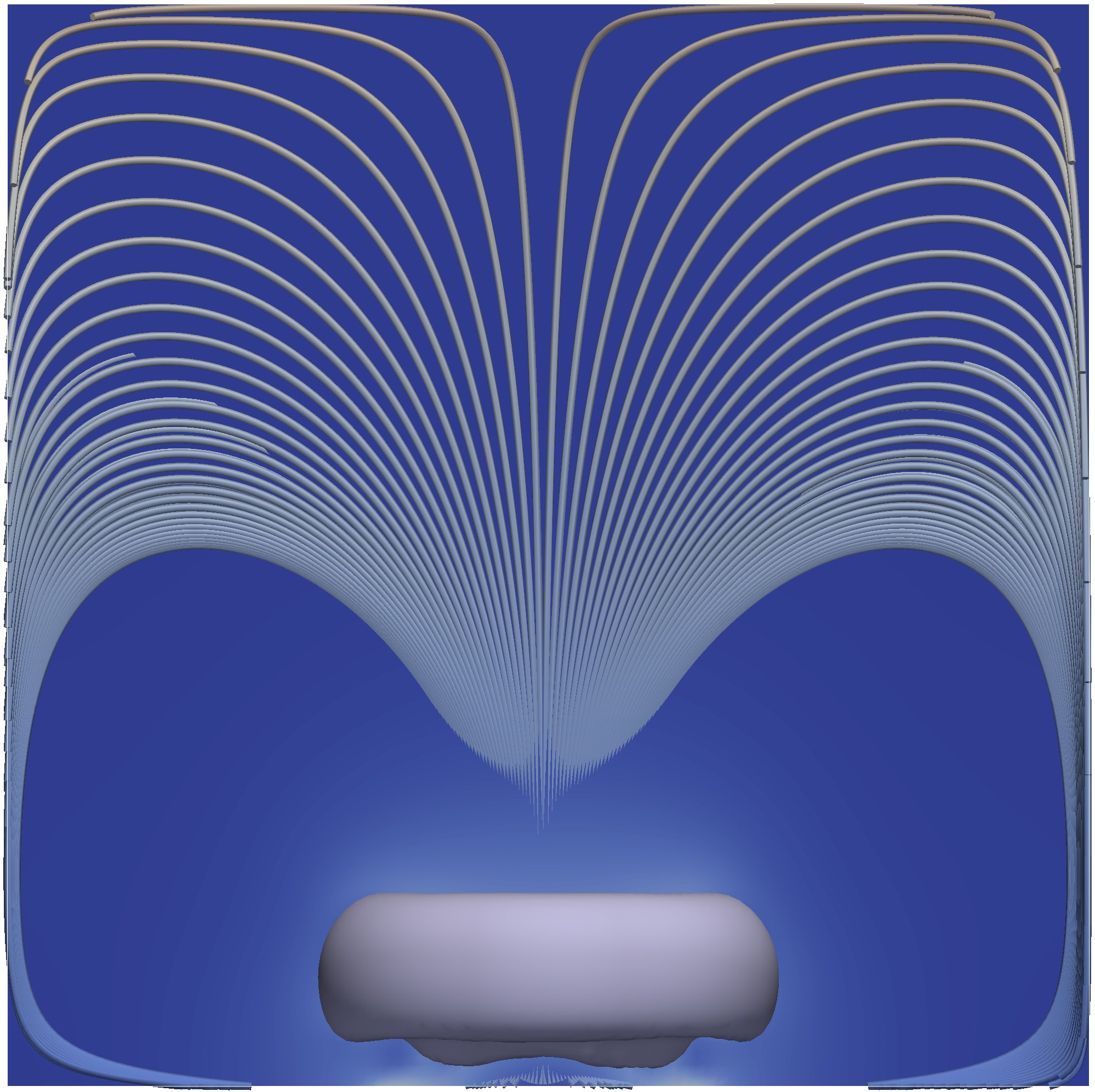}}\scalebox{0.055}{\includegraphics{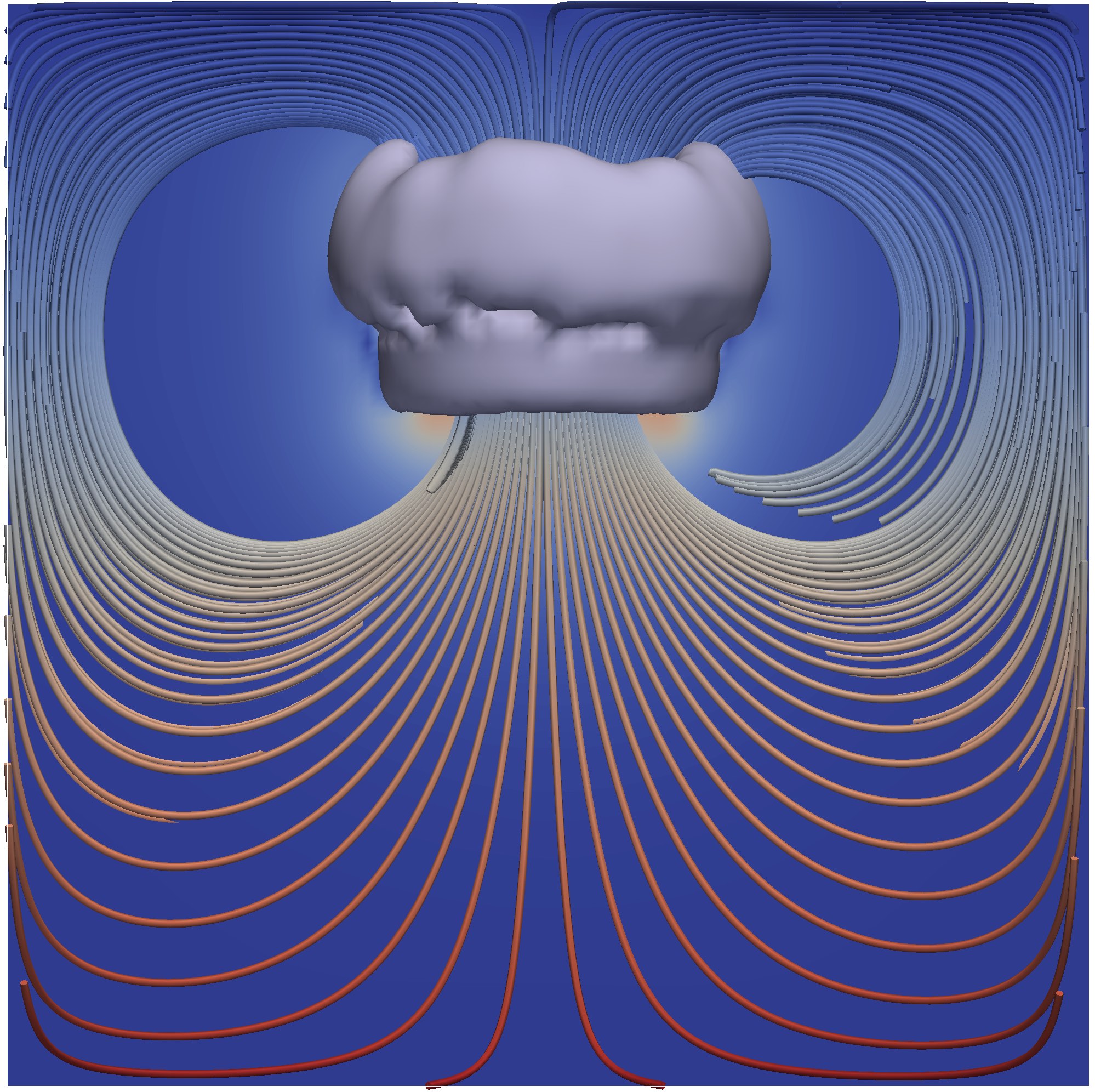}}}
\end{center}\vspace*{-1em}
\caption{\label{fig:RalizationVisualizationBubbleKL3d}Visualizations
of flow field solutions of the bubble flow with random volume force and $\pdim=3$ at $T=0.2$ seconds for stochastic
parameters $\vec{y}_{1}\approx (0.0,-0.58,-1.04)^\top$,
$\vec{y}_{12}\approx(-1.08,-1.22,-0.07)^\top$,
$\vec{y}_{19}\approx(0.97,-0.32,1.45)^\top$ and
$\vec{y}_{411}\approx(1.21,-0.93,-0.72)^\top$ (from left to right).}
\end{figure}
This two-phase flow problem is now discretized in
space by a mesh of $\Ngrid=100^3$ grid points. Time discretization is done with
a second-order Adams-Bashforth method and adaptive time step size control. The \KL expansion is truncated after $\pdim=3$ terms \reviewerone{(mainly for performance reasons)}. The
correlation length is set to $L_c=2.0$.
The norm in the radial basis function
construction is $\|\cdot\|:=\kernelscale\|\cdot\|_2$, with $\kernelscale=0.1$ for the
Gaussian kernel and all Wendland kernels and $\kernelscale=1.0$ for the Mat\'ern
kernel. 
The reference moment as in \eqref{reffirstmom} is computed by a Gaussian kernel with
$N_{\max}=1024$ collocation points.
Quadrature follows the default of tensor-product quadrature with an
approximation level of $l_q=7$.

Figure~\ref{fig:RalizationVisualizationBubbleKL3d} displays visualizations of four
solution realizations at $T=0.2$ seconds. The bubble is shown by extracting the
iso-surface of the zero level-set function.
Furthermore, a slice of the velocity field with coloring by the magnitude and
velocity field streamlines is given.

\begin{figure}
\begin{center}
\begin{minipage}[t]{0.5\textwidth}
\raisebox{2.8em}{\scalebox{0.075}{\includegraphics{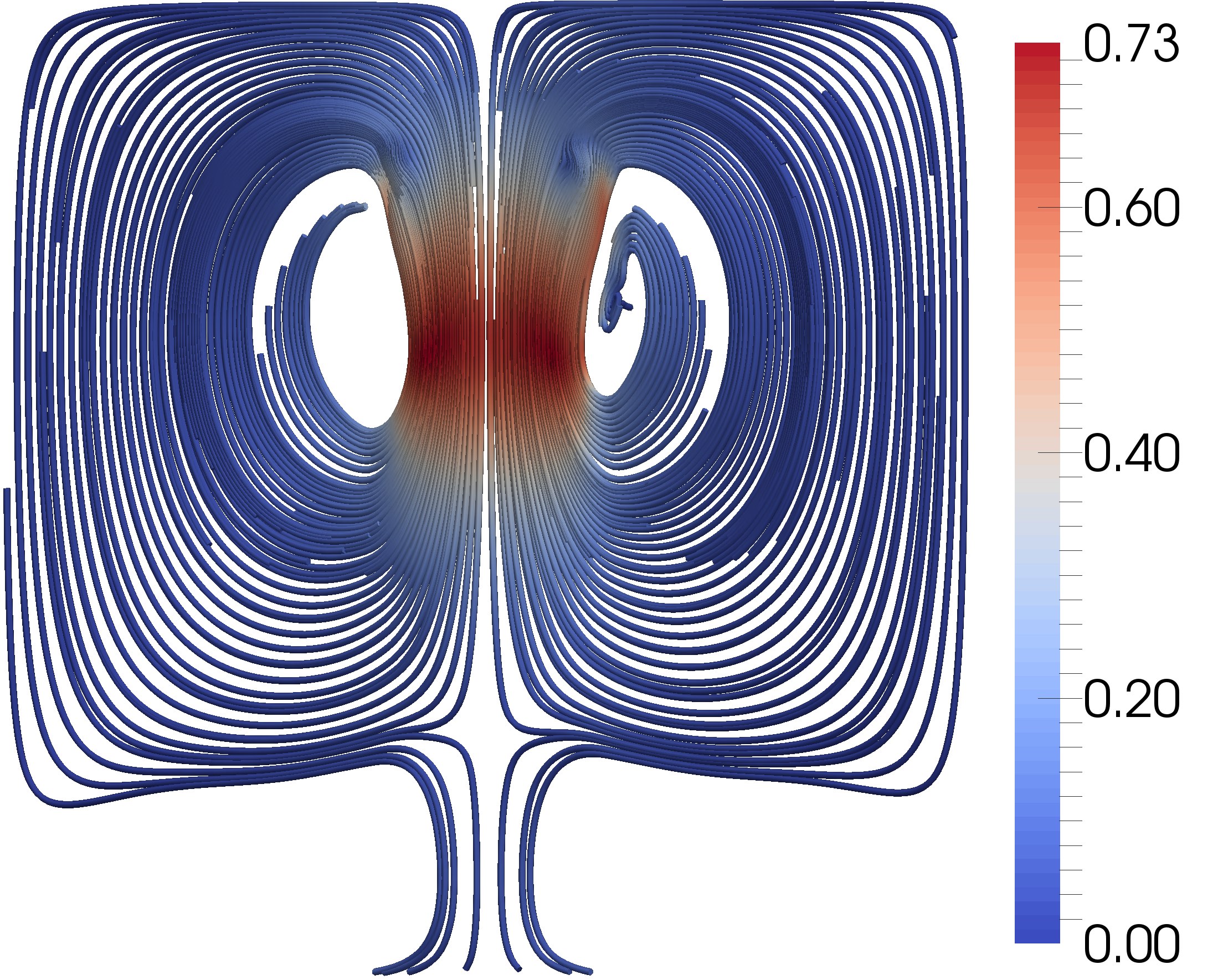}}}\\
\begin{center}
\vspace*{-4em}
\hspace*{0em}mean solution $\E{\vec{u}^{h}}(T)$
\end{center}
\end{minipage}
\scalebox{0.75}{\begin{tikzpicture}
\begin{loglogaxis}[
	height=7.8cm,
    width=7.8cm,
	xlabel=\# collocation points,
	ylabel= absolute error wrt.~reference sol.,
	legend style={legend pos=south west,font=\tiny}
	]
\addplot
table[x=collocationpoints,y=Gaussian3d]{error_mean_bubble_kl_3d_kernel_comparison.csv};
\addplot 
table[x=collocationpoints,y=Wendland3d0]{error_mean_bubble_kl_3d_kernel_comparison.csv};
\addplot 
table[x=collocationpoints,y=Wendland3d1]{error_mean_bubble_kl_3d_kernel_comparison.csv};
\addplot
table[x=collocationpoints,y=Wendland3d2]{error_mean_bubble_kl_3d_kernel_comparison.csv};
\addplot
table[x=collocationpoints,y=Wendland3d3]{error_mean_bubble_kl_3d_kernel_comparison.csv};
\addplot 
table[x=collocationpoints,y=Matern3d]{error_mean_bubble_kl_3d_kernel_comparison.csv};
\legend{Gauss. $\kernel_\epsilon$, Wend. $\kernel_{3,0}$,
Wend. $\kernel_{3,1}$, Wend. $\kernel_{3,2}$,
Wend. $\kernel_{3,3}$, \parbox{1.6cm}{Mat\'ern
$\kernel_{\frac{3+3}{2}}$}}
\addplot[black] coordinates {
                (16, 4.0e-3)
               (128, 8.40896416e-4)
        };
\draw (axis cs:128,8.40896416e-4) node[below]{{\small $0.75$}};
\end{loglogaxis}
\end{tikzpicture}}
\end{center}\vspace*{-1em}
\caption{\label{fig:ResultsErrorMeanBubbleKL3dKernelOverview} Streamline
slice visualization of the mean velocity field with color-coded velocity magnitude in the bubble flow problem with random volume
forces (left) and error convergence results for the first component of the mean velocity field $\E{u^{h}_1}(T)$ (right).}
\end{figure}

On the left-hand side in
Figure~\ref{fig:ResultsErrorMeanBubbleKL3dKernelOverview}, a visualization of the mean velocity field $\E{\vec{u}^{h}}$ by means of a slice of streamlines through all three mean velocity field components is shown.
The diagram in Figure~\ref{fig:ResultsErrorMeanBubbleKL3dKernelOverview} (right) displays the convergence of the approximation of $\E{u^{h}_1}(T)$ for different kernel functions with respect to the reference solution. We observe that rates in the range of $0.75$ are achieved by the Wendland and Mat\'ern kernels. 
There is a slight reduction of the convergence rate for smoother Wendland
kernels. However, the convergence
results obtained by the Gaussian kernel clearly show higher-order and maybe even an exponential convergence behavior might be anticipated here. \reviewertwo{This suggests a smooth dependence of the quantity of interest on the random input.}

\subsection{Rising bubble flow in a random situation}
\label{sec:ModelsStochasticHomogenizationForRisingBubbles}
Our second application for the two-phase Na\-vier-Stokes equations uses a similar setup as the first flow example. We again consider a rising air bubble in some liquid. Now, however, the density, viscosity and
initial bubble position are under stochastic influence, whereas the volume forces are
deterministically given. 

Again Figure~\ref{fig:ModelRisingBubbleSetup} outlines the basic setup
of this two-phase flow problem. The domain is given as $\bar{\dom}=[0,0.2]^3$
and we have $T=0.35$.The volume force $\vec{g}$ is given as standard gravity.
The initial conditions are also set to
$\vec{u}_{0_i}=\left(0,0,0\right)^\top$. Furthermore, the surface tension
coefficient is $\sigma=0.0728$. The boundary conditions are the same as in the
previous test case.
The gas phase $\dom_{2}$ in once more modeled as air, i.e., 
the density is $\rho_2=1$, and the viscosity is $\mu_2=1.72\cdot 10^{-5}$. The remaining quantities are assumed to be random. The missing random density $\rho_{1}$ and viscosity $\mu_{1}$ will be given in \eqref{randompara}.
The initial position of
the air bubble shall be a random quantity, see \eqref{randompos}. This means that the phase-wise sub-domains
$\dom^{h}_1,\dom^{h}_2$ (or the level-set function describing them) are stochastic processes themselves. 
Therefore, these domains are given for time $t$ by\vspace*{-0.5em}
$\dom^{h}_i(\vec{y},t)$ with
\begin{linenomath*}
\reviewerone{\begin{equation}
\dom^{h}_i(\vec{y},t) = \Phi[(\dom^{h}_i)^0(\vec{y})](t),\quad (\dom^{h}_i)^0(\vec{y}) :=
\dom^{h}_i(\vec{y},0),\quad y\in\paradom\,,
\end{equation}}
\end{linenomath*}
where $\Phi$ again describes the transformation of the initial domains
$(\dom^{h}_i)^0(\vec{y})$ under fluid flow. It deterministically depends on the
velocities $u_1,u_2$, but the initially given domains are subject to random
perturbations. We define the initial liquid and gas phase domains
$\dom_1^0(\vec{y})$ and $(\dom^{h}_2)^0(\vec{y})$ such that the gas phase domain is
a sphere of radius $0.03\,m$ around some random initial center
$\vec{x}_{init}(\vec{y})$ at the beginning, thus\vspace*{-0.5em}
\begin{linenomath*}
\reviewerone{
\begin{equation}
(\dom^{h}_1)^0(\vec{y}) := \left\{ \vec{x}\in\dom
\middle| \|\vec{x}-\vec{x}_{init}(\vec{y})\| > 0.03\right\}\,,
\end{equation}
}
\end{linenomath*}
\begin{linenomath*}
\reviewerone{\begin{equation}
(\dom^{h}_2)^0(\vec{y}) := \left\{
\vec{x}\in\dom \middle| \|\vec{x}-\vec{x}_{init}(\vec{y})\| < 0.03\right\}\,.
\end{equation}}
\end{linenomath*}
The initial free surface is $\bnd_f(\vec{y},0) =
\dom^{h}\setminus((\dom^{h}_1)^0(\vec{y})\cup (\dom^{h}_2)^0(\vec{y}))$. Moreover, the
random parameter functions
$x_1^{init}(\vec{y})$, $x_2^{init}(\vec{y})$ and $x_3^{init}(\vec{y})$ with
$\vec{x}_{init}(\vec{y})=\left(x_1^{init},x_2^{init},x_2^{init}\right)$ as
well as the material parameters for the liquid phase, $\mu_1(\vec{y})$ and
$\rho_1(\vec{y})$ are modeled by truncated \KL expansions. Truncation is done
after the first stochastic term. Overall, these functions are given as
\begin{linenomath*}
\begin{align}
&x_1^{init}(\vec{y}) =  0.1 + \frac{0.06}{\sqrt{3}} y_1^{x_1}\,, \quad 
x_2^{init}(\vec{y}) =  0.06 + \frac{0.01}{\sqrt{3}} y_1^{x_2}\, \quad
x_3^{init}(\vec{y}) = 0.1 + \frac{0.06}{\sqrt{3}} y_1^{x_3}\,,\label{randompos}\\
&\mu_1(\vec{y}) = 0.5005 + \frac{0.4995}{\sqrt{3}} y_1^{\mu_1}\,, \quad 
\rho_1(\vec{y}) = 750 + \frac{250}{\sqrt{3}} y_1^{\rho_1}\,. \label{randompara}
\end{align}
\end{linenomath*}

All parameters of this problem are now collected in the five-dimensional vector $\vec{y}$ with $\vec{y}=\left(y_1^{x_1},y_1^{x_2},y_1^{x_3},y_1^{\mu_1},y_1^{\rho_1}\right)^\top$.
We assume $\vec{y}$ to be a realization of the independent random vector $Y$ with
\begin{linenomath*}
\reviewerone{\begin{equation}
Y_1^{x_1}\sim\mathcal{U}\left(-\sqrt{3},\sqrt{3}\right), \quad
Y_1^{x_2}\sim\mathcal{U}\left(-\sqrt{3},\sqrt{3}\right), \quad
Y_1^{x_3}\sim\mathcal{U}\left(-\sqrt{3},\sqrt{3}\right)\,,
\end{equation}}
\reviewerone{\begin{equation}
Y_1^{\mu_1}\sim\mathcal{U}\left(-\sqrt{3},\sqrt{3}\right), \quad
Y_1^{\rho_1}\sim\mathcal{U}\left(-\sqrt{3},\sqrt{3}\right)\,.
\end{equation}}
\end{linenomath*}
Consequently, the density function becomes
$\rho(\vec{y})=\left(\frac{1}{2 \sqrt{3}}\right)^5$. Here, we again approximate
the first stochastic moment $\E{u^{h}_1}(T)$.

This flow problem is now discretized using a uniform grid with $\Ngrid=100^3$ grid points, a finite-difference discretization in space and a
second-order Adams-Bashforth method with adaptive time-stepping in time. 
\begin{figure}
\begin{center}\scalebox{1.06}{\scalebox{0.055}{\includegraphics{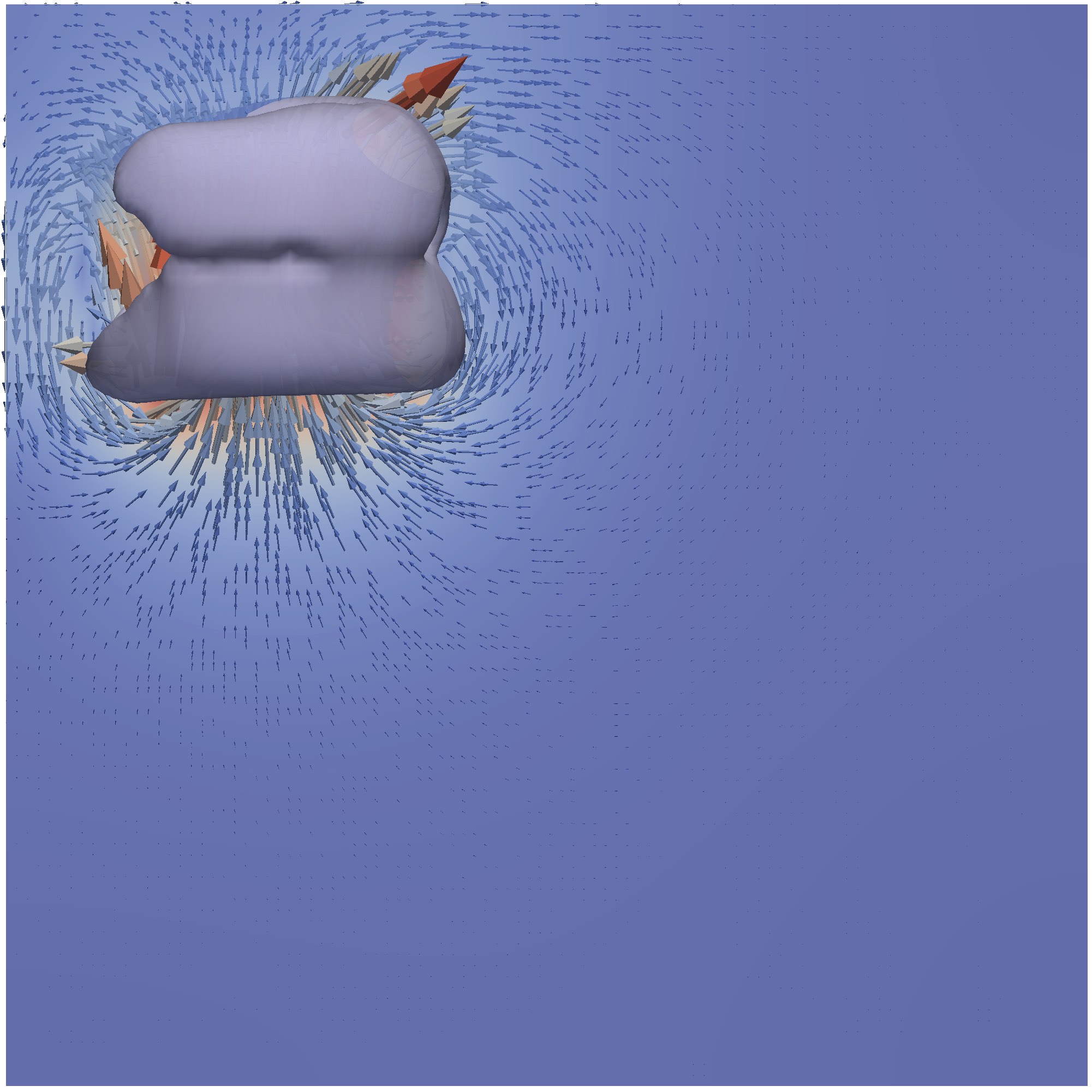}}\scalebox{0.055}{\includegraphics{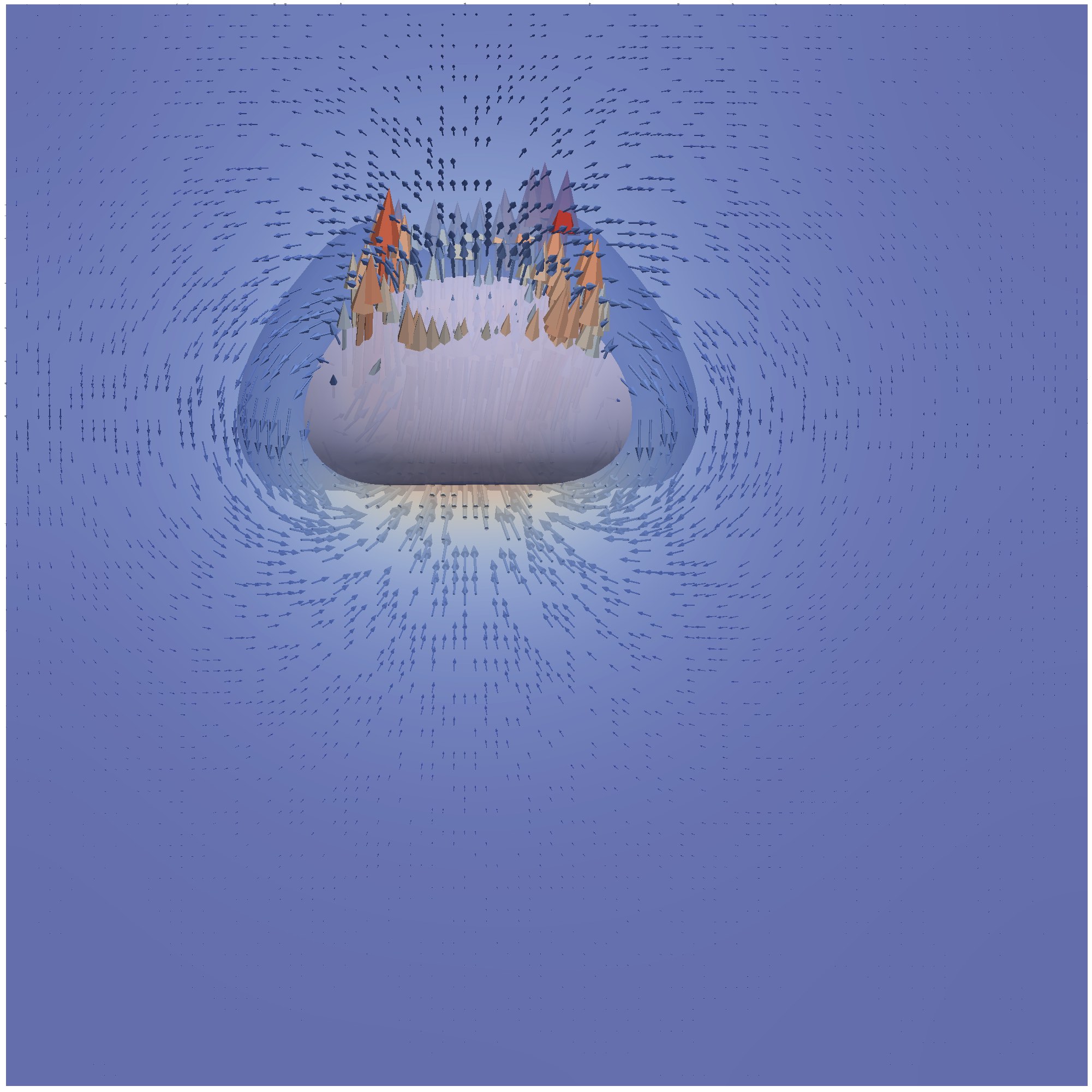}}\scalebox{0.055}{\includegraphics{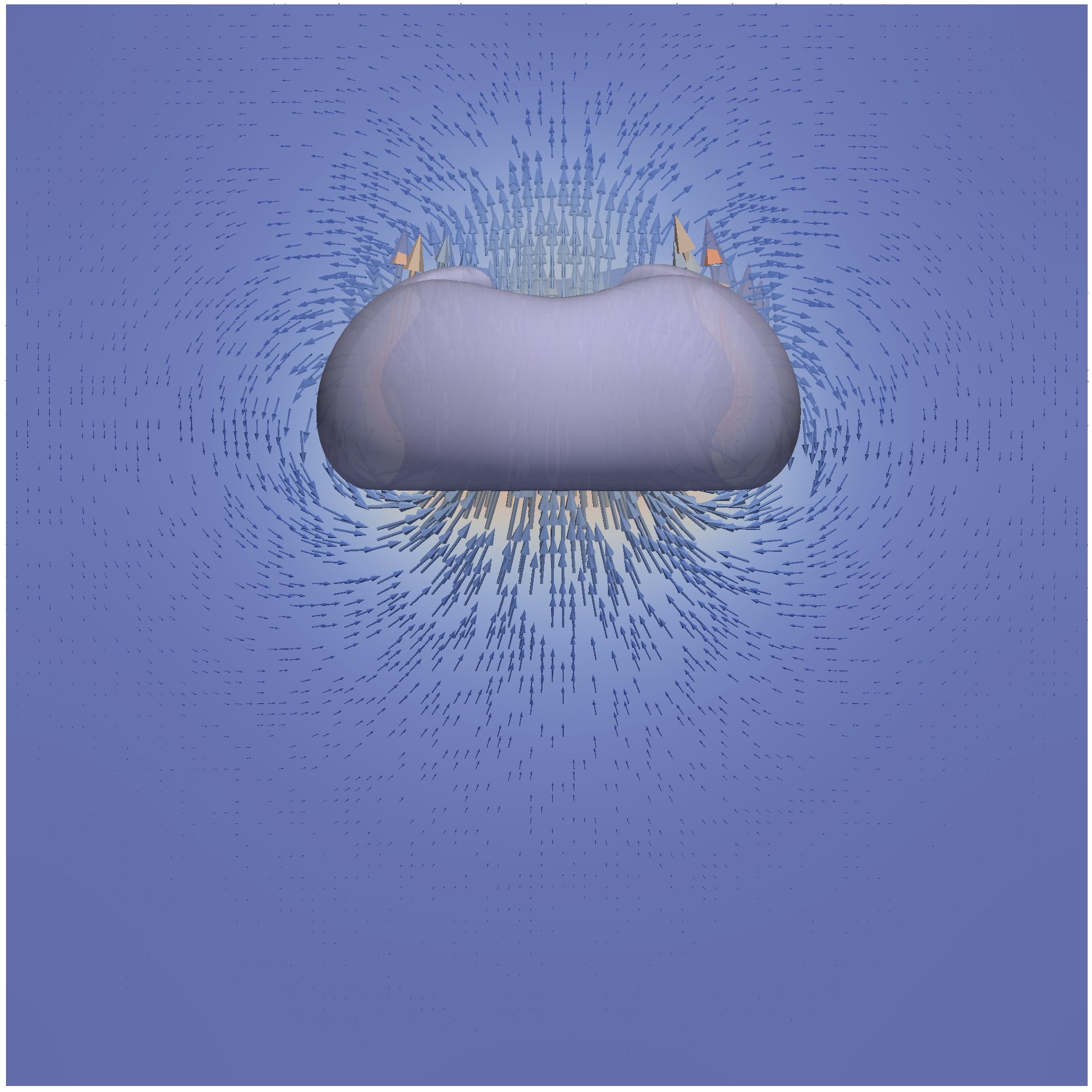}}\scalebox{0.055}{\includegraphics{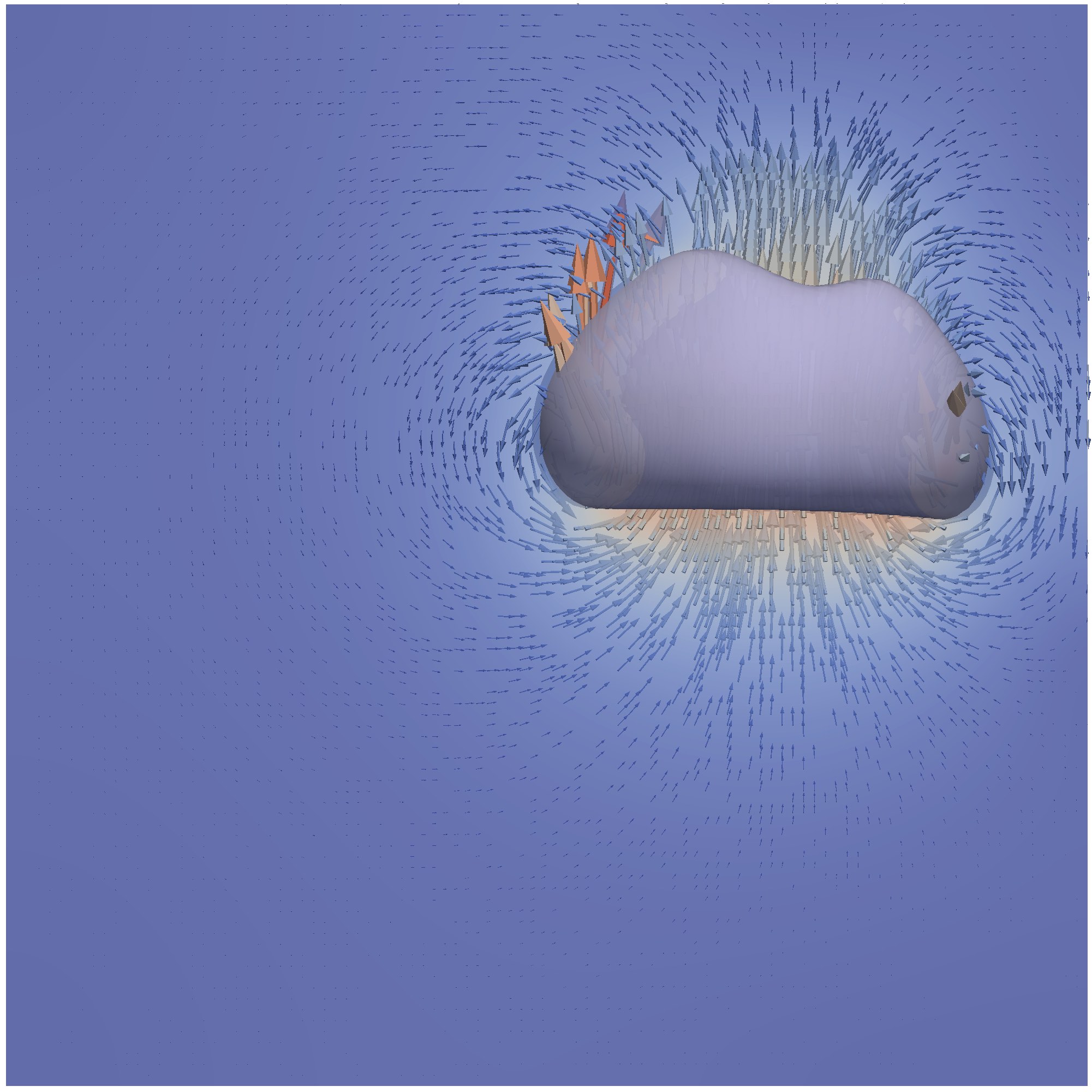}}}
\end{center}\vspace*{-1em}
\caption{\label{fig:RealizationsVisualizationBubbleHomogenization}Flow
field and bubble visualization of solution realizations of the rising bubble problem with stochastic parameters
$\vec{y}_{8}^\prime\approx (0.05,0.07,0.12,0.02,863.64)^\top$,
$\vec{y}_{6}^\prime\approx(0.09,0.05,0.07,0.09,772.73)^\top$,
$\vec{y}_{33}^\prime\approx(0.1,0.05,0.12,0.08,512.4)^\top$
and 
$\vec{y}_{27}^\prime\approx(0.14,0.05,0.09,0.09,735.54)^\top$ (from
left to right).}
\end{figure}
We apply the stochastic
approximation in the image of the five random input \reviewerone{variables} (see also \eqref{randompos} and \eqref{randompara}), thus a stochastic
space $\paradom^\prime=[0.04,0.16]\times[0.05,0.07]\times [0.04,0.16]\times
[0.001,0.1]\times [500,1000]$ is considered. To approximate this space
uniformly, the norm
\begin{linenomath*}
\reviewerone{\begin{equation}
\left\|\vec{y}^\prime\right\| := \kernelscale\left\|(5 x_1^{init},
5 x_2^{init}, 5 x_3^{init}, \mu_1, 10^3 \rho_1)\right\|_2
\end{equation}}
\end{linenomath*}
is used in the construction of the radial basis functions. Approximation by
Gaussian kernels is done with $\kernelscale=0.1$, while the other kernels have
$\kernelscale=1.0$. All other approximation parameters are set
as mentioned in Section~\ref{sec:NumResBenchmarkSetup}. 

Figure~\ref{fig:RealizationsVisualizationBubbleHomogenization} shows four flow
field realizations at $T=0.2\, s$. 
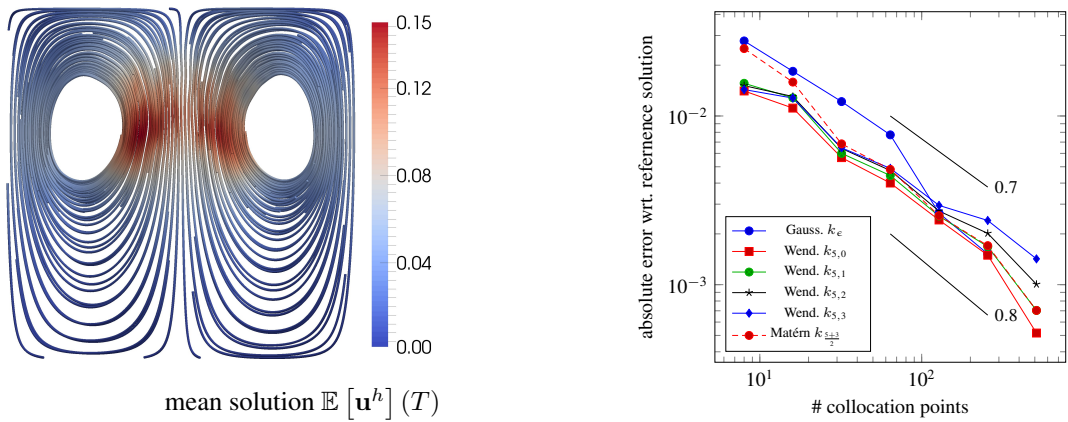
\begin{figure}
\begin{center}
\begin{minipage}[t]{0.5\textwidth}
\raisebox{2.4em}{\scalebox{0.075}{\includegraphics{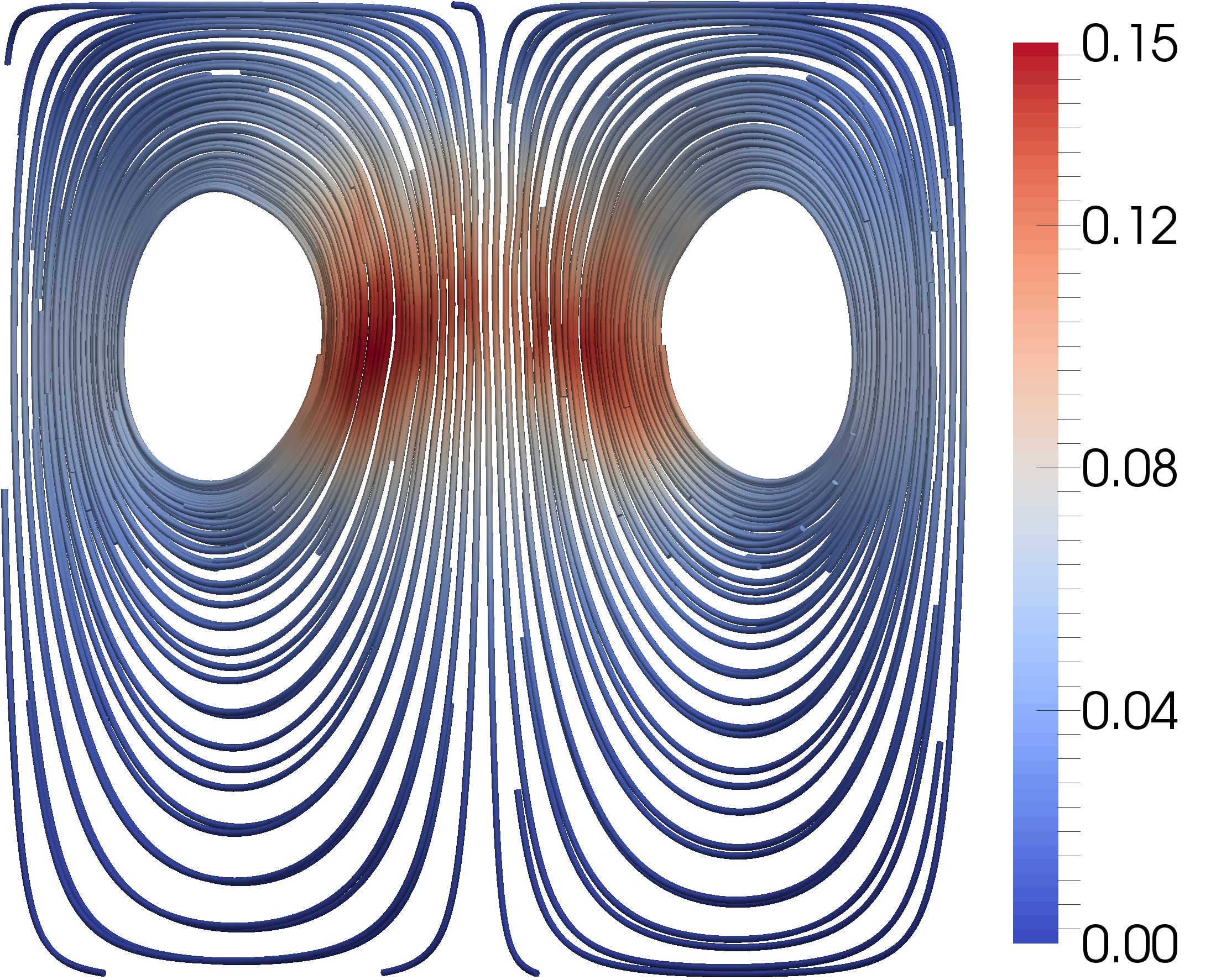}}}\\
\begin{center}
\vspace*{-3.7em}\hspace*{-1em}mean solution $\E{\vec{u}^{h}}(T)$
\end{center}
\end{minipage}
\scalebox{0.75}{\begin{tikzpicture}
\begin{loglogaxis}[
	height= 7.8cm,
    width= 7.8cm,
	xlabel=\# collocation points,
	ylabel= absolute error wrt.~refernence solution,
	legend style={legend pos=south west,font=\tiny}
	]
\addplot
table[x=collocationpoints,y=Gaussian3d]{error_mean_bubble_homogenization_kernel_comparison.csv};
\addplot 
table[x=collocationpoints,y=Wendland3d0]{error_mean_bubble_homogenization_kernel_comparison.csv};
\addplot 
table[x=collocationpoints,y=Wendland3d1]{error_mean_bubble_homogenization_kernel_comparison.csv};
\addplot
table[x=collocationpoints,y=Wendland3d2]{error_mean_bubble_homogenization_kernel_comparison.csv};
\addplot
table[x=collocationpoints,y=Wendland3d3]{error_mean_bubble_homogenization_kernel_comparison.csv};
\addplot 
table[x=collocationpoints,y=Matern3d]{error_mean_bubble_homogenization_kernel_comparison.csv};
\legend{Gauss. $\kernel_\epsilon$, Wend. $\kernel_{5,0}$,
Wend. $\kernel_{5,1}$, Wend. $\kernel_{5,2}$,
Wend. $\kernel_{5,3}$,
\parbox{1.6cm}{Mat\'ern $\kernel_{\frac{5+3}{2}}$}}
\addplot[black] coordinates
{ (64, 2.0e-3)
               (256, 0.659753955e-3)
        };
\draw (axis cs:256,0.659753955e-3) node[right]{{\small $0.8$}};
\addplot[black] coordinates {
                (64, 1.0e-2)
               (256, 3.78929142e-3)
        };
\draw (axis cs:256,3.78929142e-3) node[right]{{\small $0.7$}};
\end{loglogaxis}
\end{tikzpicture}}
\end{center}\vspace*{-1em}
\caption{\label{fig:ResultsErrorMeanBubbleHomogenizationKernelOverview}
\reviewertwo{Left: Streamline slice visualization of the mean velocity field (color-coded by velocity magnitude)
in the bubble flow stochastic homogenization problem. Right: Error convergence results for the
approximated mean of the first component of the velocity field $\E{U_1}(T)$ of the same problem.
Approximations using different kernel functions are compared.}}
\end{figure}
Approximation results for the mean of the velocity field at $T=0.2\, s$ are
given in~Figure~\ref{fig:ResultsErrorMeanBubbleHomogenizationKernelOverview}. 
One the left, a streamline slice visualization of the reference
solution is shown, which is approximated by the Gaussian kernel with $N_{\max}
=512$
collocation points. On the right-hand side, error convergence in the mean $\E{u^{h}_1}(T)$ of the
first component of the velocity field is given for different kernel
functions. We observe that rates in the range of $0.7$--$0.8$ are achieved by Wendland, Mat\'ern and Gauss kernels. The use of higher-order Wendland kernels does not result in any improvement. The highest convergence rates are achieved by the Gaussian and the Mat\'ern kernel with
about $0.8$.

\subsection{Comparison to sparse spectral tensor-product approximations}\vspace*{-0.5em}
\label{sec:NumResDakotaComparison}
In this section, our proposed kernel-based stochastic collocation method is
compared to sparse grid approximations. 
To this end, we employ on the one hand our kernel-based method with the default parameter chosen
as in Section~\ref{sec:NumResBenchmarkSetup},
and on the other hand the sparse grid method within the \textit{Dakota} framework \cite{Adams2009a}.
Dakota is a parallel software suite developed by the Sandia
National Laboratories.
It allows e.g.~to apply optimization and uncertainty quantification for
black-box solvers.
Dakota features, among others, 
sparse grid stochastic collocation.

\begin{figure}
\begin{center}
\scalebox{0.75}{\begin{tikzpicture}
\begin{loglogaxis}[
	height=7.8cm,
    width=7.8cm,
    xmax=650,
	xlabel=\# collocation points,
	ylabel=relative error wrt.~reference sol., 
	legend style={legend pos=south west, font=\tiny}
	]
\addplot 
table[x=collocationpoints,y=Gaussian3d]{error_mean_bubble_kl_general_comparison.csv};
\addplot 
table[x=collocationpoints,y=Matern3d]{error_mean_bubble_kl_general_comparison.csv};
\addplot
table[x=collocationpoints,y=StochCollocSparseGrid3d]{error_mean_bubble_kl_general_comparison.csv};
\legend{
Gauss. $\kernel_\epsilon$,
Mat\'ern $\kernel_{\frac{3+3}{2}}$,
\reviewertwo{sparse grid stoch.~colloc.}}
\end{loglogaxis}
\end{tikzpicture}}
\end{center}\vspace*{-1em}
\caption{\label{fig:ResultsErrorMeanBubbleKLGeneralComparison}\reviewertwo{Error convergence results of the approximation of the mean
bubble center position $\E{c^{h}_2}(T)$ in the large-scale two-phase flow problem with rising bubbles. Kernel-based stochastic collocation (with Gaussian and Mat\'ern kernel) and sparse grid-based stochastic collocation are compared.}}
\end{figure}
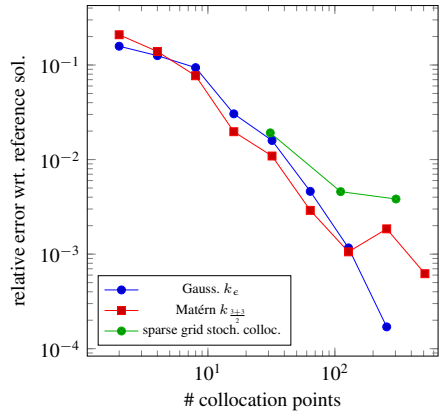

Classical stochastic collocation \reviewerone{as provided by Dakota} uses univariate Lagrange polynomials as
Lagrange basis functions \cite{Babuska2010}. 
For multi-variate interpolation, i.e., higher dimensions in stochastic space, 
a Smolyak sparse grid construction \reviewerone{as in \cite{Babuska2010}}
is used. \reviewerone{More details on the Dakota-based implementation of the Smolyak sparse grid construction are given in \cite{Adams2009a}.  
The employed parameters for sparse grid stochastic collocation in Dakota 
are depicted in Configuration~\ref{lst:NumResSCDakotaConfig} in~\ref{appendix1}.
We will employ the Dakota implementation for the numerical comparison in Figure~\ref{fig:ResultsErrorMeanBubbleKLGeneralComparison}.
}

Note here that, we do not consider \textit{an\-iso\-tro\-pic} sparse grid constructions
\cite{Nobile2008a}, because the RBF kernel-based stochastic
collocation is employed without any directional preference anyway.
However, the topic of dimension-wise weighting is future work and has been already partially considered in \cite{Zaspel2015}.

For reason of simplicity, in contrast to the previous paragraphs, the subsequent convergence
study compares kernel-based results and the Dakota-based sparse grid result using only the
single-valued quantity of interest,\linebreak $\E{c^{h}_2}(T)$. The single deterministic two--phase flow problems are discretized and solved as before by using\linebreak NaSt3DGPF. In the Dakota calculations, NaSt3DGPF
is directly called from the Dakota control program. 

We come back to the problem of a rising gas bubble in water which is subject
to a random volume force,
cf.~Section~\ref{sec:ModelsKLbasedRandomVolumeForceInBubbleFlow}.
Again, the random input is modeled by a \KL expansion which is truncated after the third
term, thus we have $\pdim=3$ stochastic dimensions and a correlation length of
$L_c=2.0$ is used. The quantity of interest is now
the second component $c^{h}_2$ of the bubble's center position at physical time $T=0.2$ seconds.
The reference mean $\E{c^{H}_2}(T)$ is approximated by a Gaussian
kernel with $N_{\max}=512$ collocation points. 
In case of kernel-based approximation, the standard RBF norm $\|\cdot\|:=\kernelscale
\|\cdot\|_2$ with a modified scaling of $\kernelscale=0.1$ is employed. 
Results computed with the
Gaussian kernel are regularized with $\epsilon_{reg}=10^{-6}$ while Mat\'ern
kernel results are regularized with $\epsilon_{reg}=10^{-5}$. Quadrature is the full tensor-product quadrature on level $l_q=9$. 

In Figure~\ref{fig:ResultsErrorMeanBubbleKLGeneralComparison}, convergence
results are given for a mean approximation by Gaussian kernels, Mat\'ern
kernels and stochastic collocation by sparse grids (SG). We observe that the kernel-based methods show better convergence than the sparse grid method and a
higher-order algebraic rate is achieved. Moreover, we even observe a spectral convergence rate for the Gaussian kernel. 
Furthermore, the errors of the RBF-based stochastic collocation method are always below the results of the sparse grid stochastic
collocation method. 

Note at this point that a computationally challenging problem is
considered and to achieve an error with almost the size of $10^{-4}$ with only 256
simulations is quite remarkable. It is well-known that sto\-chastic collocation on sparse tensor product
constructions shows asymptotically exponential convergence rates on sufficiently
smooth problems, see \cite{Babuska2010}. Nevertheless, in many large-scale uncertainty quantification application problems,
the quantity of interest has either limited smoothness with respect to the
random input and/or bad pre-asymptotic error behavior is dominating the convergence to a significant extend. This is where
kernel-based stochastic collocation has its main advantage.

\section{Conclusions}
\label{sec:conclusions}
In this work, the solution of random two-phase Navier-Stokes problems by means of the 
kernel-based stochastic collocation has been considered. Here, the first stochastic moment of
solutions of two--phase flow problems was computed by applying a new \textit{non-intrusive method}
to the existing flow solver NaSt3DGPF which is based on a kernel-approximation. 
For any method, given a fixed
target error tolerance, it is indispensable to keep the number of
discretization points in stochastic space as low as possible. A way to
overcome this issue is to introduce an approximation method in stochastic
space, which has high convergence order with a very small
pre-asymptotic error. This has been achieved in this work by the introduction
of the RBF kernel-based stochastic collocation method. 
Numerical results were given that underline the good properties of the
kernel-based method. For the 
random two-phase Navier-Stokes equations, algebraic convergence rates were shown and a small error in the pre-asymptotic regime was always present.
Thus, our kernel--based stochastic colocation approach outperformed well-known
established methods such as Monte Carlo, or sparse
spectral tensor-product stochastic collocation in this situation.

\appendix
\section{Parameters for sparse grid stochastic collocation}\label{appendix1}
\lstset{caption={Standard parameters for sparse grid stochastic collocation in Dakota},label=lst:NumResSCDakotaConfig}
\begin{lstlisting}
method
  stoch_collocation
    sparse_grid_level = <level>
    dimension_preference = <N_KL>*1
    samples = 10000 seed = 12347 rng rnum2
    output silent
\end{lstlisting}

\acknowledgements
Major parts of the numerical results were computed on a GPU cluster of the
institute SCAI which is part of the Fraunhofer Society. 
This support is gratefully acknowledged.
\bibliographystyle{plain}

\end{document}